\documentclass{article}

\usepackage{amsmath,amssymb,amsthm,mathrsfs,graphicx}
\usepackage[titletoc, title]{appendix}
\usepackage{mathrsfs}
\usepackage{paralist}
\usepackage{graphics} 
\usepackage{epsfig} 
\usepackage{multirow}
\usepackage{graphicx}
\usepackage{epstopdf}
\usepackage{syntonly,indentfirst,mathrsfs,latexsym,float}
\usepackage{bm,amsmath,amsbsy,amsthm,amssymb,amsfonts}
\usepackage[colorlinks,linkcolor=red]{hyperref}
\usepackage[capitalise, nosort]{cleveref}
\usepackage{array}
\usepackage{color}
\usepackage{multirow}
\usepackage{mathrsfs}
\usepackage{type1cm}
\usepackage{wasysym}

\numberwithin{equation}{section}

\usepackage{indentfirst} 
\usepackage{stmaryrd}
\usepackage{subfigure}
\usepackage[utf8]{inputenc}
\usepackage{pgf,tikz}
\usepackage{threeparttable}
\usetikzlibrary{arrows}

\newtheorem{thm}{Theorem}[section]

\newtheorem{Def}[thm]{Definition}
\newtheorem{Lem}[thm]{Lemma}
\newtheorem{Rem}{Remark}[section]
\newtheorem{Exa}{Example}[section]

\usepackage{authblk}
\usepackage{enumerate}

\usepackage{lineno}    

\title{{Semi-discrete and fully discrete weak Galerkin finite element methods  for a quasistatic Maxwell viscoelastic model }\thanks{This work was supported  by the National Natural Science Foundation of China (12171340)}}

\author[a,b]{Jihong Xiao\footnote{Email: xiaojh2752@163.com}}
\author[a]{Zimo Zhu\footnote{ Email: zzm@stu.scu.edu.cn}}
\author[a]{Xiaoping Xie \thanks{Corresponding author. Email: xpxie@scu.edu.cn}}
\affil[a]{School of Mathematics, Sichuan University, Chengdu 610064,China}
\affil[b]{Mathematics department of Jinjiang college, Sichuan University, Pengshan 620860,China}

\date{}

\begin{document}
\maketitle
	
\begin{abstract} \large	\noindent{}$\mathbf{Abstract}$: This paper  considers   weak Galerkin  finite element approximations for a quasistatic Maxwell viscoelastic model. The spatial discretization uses piecewise polynomials of degree  $k \ (k\geq 1)$ for the  stress approximation,  degree $k+1$ for the velocity approximation, and degree $k$ for the numerical trace of velocity  on the inter-element boundaries. The temporal discretization in the fully discrete method adopts a backward Euler difference scheme. We show  the existence and uniqueness of the semi-discrete and fully discrete solutions, and derive optimal a priori error estimates.  Numerical examples are provided to support the  theoretical analysis. \\[2mm]	

\end{abstract}
	
	\noindent{} $\mathbf{Keywords}$:  quasistatic Maxwell viscoelastic model; Weak Galerkin method; Semi-discrete scheme; Fully discrete scheme; Error estimate. 
	\\[2mm]
	\noindent{} $\mathbf{Mathematics\,\, subject\, classifications(2020): 35Q74, 65M12,  65M60.}$ 
	
	\section{Introduction} 
Let $\Omega\subset \mathbb{R}^d(d=2,3)$ be a convex polyhedral domain with boundary $\partial \Omega$, and $T$ be a positive constant. 
We consider  the following quasistatic Maxwell viscoelastic model:
 \begin{subequations}	\label{modelproblem}
\begin{align}
  -\mathbf{div}\bm {\sigma} =&\bm{f} , &  &  (x,t)\in\Omega\times[0,T] ,  \label{mainequation}
 \\ \bm\sigma +\bm\sigma_t  =&   \mathbb{C}\bm{\varepsilon(u_t)}, & & (x,t)\in\Omega\times[0,T],   \label{modelproblem2}
\\ \bm{u}  = &0, & &(x,t)\in\partial\Omega\times[0,T], 
\\ \bm{u}(x,0)=&\phi_0(x), & &x\in \Omega, 
\\ \bm\sigma(x,0)=&\psi_0(x), & &x\in \Omega.
\end{align}
\end{subequations}
Here $\bm{u}\in \mathbb{R}^d$ is the displacement field, $\bm\sigma=(\sigma_{ij})_{d\times d} $ the symmetric stress tensor, $\bm{\varepsilon(u)} = \left(\bm{\nabla u + (\nabla u)^T}\right)/2$ the strain tensor,   $f$   the body force,   $\phi_0(x)$ and $\psi_0(x)$ are initial data,   $g_t:=\partial g/\partial t $  for any function $g(\boldsymbol{x},t)$, and $\mathbb{C}$ denotes an elastic module tensor satisfying
\begin{align}   \label{assumed_condition}
    0 < M_0\bm{\tau}:\bm{\tau}\leq \mathbb{C}^{-1}\bm\tau:\bm\tau
   \leq M_1 \bm\tau:\bm\tau , \forall \hspace{1mm}symmetric\hspace{1mm} tensor\hspace{1mm} \bm{\tau}=(\tau_{ij})_{d\times d}  \ a.e. \hspace{2mm} x\in \Omega, 
\end{align}
where $M_0$ and $ M_1$ are two positive constants, and $\bm{\iota:\tau}:=\sum\limits_{i=1}^d\sum\limits_{j=1}^d {\iota_{ij}\tau_{ij}}$ for   $\bm{\iota},{\bm\tau}\in \mathbb{R}^{d\times d}$.
Note that for an  isotropic elastic medium we have
$$\mathbb{C}\bm{\varepsilon(u_t)} = 2\mu\bm{\varepsilon(u_t)} +\lambda (\nabla\cdot \bm{u_t}) I,$$
where $\mu $ and $\lambda$ are   Lam$\acute{e}$ constants, and $I$ the identity matrix.  

In material science and continuum mechanics, viscoelasticity is the property of materials that exhibit both viscous and elastic characteristic when undergoing deformation.    The Maxwell model,   characterized by the governing constitutive relation \eqref{modelproblem2},  is one of   classical models of viscoelasticity; see, e.g. 
\cite{Bland1960,   Christensen,
Dill2007, Drozdov1998,  Fung1966International,  Gold1988,  Gurtin1962,  Marques2012,  Salencon2016, Schapery2000} for some related works on  the development and applications of viscoelasticity theory. These models, including the Kelvin-Voigt model and the Zener model,  are represented by different combinations of  purely elastic springs, which obey   Hooke's law,  and purely viscous  dashpots, which  obey Newton law.   The Maxwell model consists of a    spring and  a dashpot connected  in series. We note that the general constitutive law of viscoelasticity can be described in a unified framework   by using convolution integrals in time with some kernels \cite{Christensen, Fung1966International, Salencon2016}.

In \cite{Carcione1988A, Carcione1988B} Carcione et al.  gave the first numerical simulation of wave propagation in viscoelastic materials, and 
 introduced  memory variables to avoid the computation of convolution integrals in the constitutive relation.
 Janovsky et al.  \cite{Janovsky1995} applied continuous/discontinuous Galerkin finite element  methods to discretize a linear viscoelasticity model involving the hereditary constitutive relations for compressible solids. 
Ha et al.\cite{Ha2002} proposed a nonconforming finite element method for a viscoelastic complex model in the space frequency domain.
B\'ecache el al. \cite{Becache2002} presented   a family of mass lumped mixed finite element methods,   together with a leap-frog scheme  in the time discretization, for the Zener model. In \cite{Riviere2006, Riviere2003, Riviere2007} Rivi\`{e}re et al. analyzed   discontinuous Galerkin finite element discretizations of the quasistatic linear viscoelasticity and linear/nonlinear diffusion viscoelastic models, where   a Crank-Nicolson temporal scheme is used in the full discretization.
Rognes and Winther  \cite{Rognes2010}
considered mixed finite element approximations with weak symmetric stresses for the quasistatic Maxwell and Kelvin-Voigt models, where the temporal discretization uses a second backward difference scheme. In \cite{Shidongyang2011} Shi and Zhang applied the standard $p$-order rectangular finite elements   to solve a kind of nonlinear viscoelastic wave equations with nonlinear boundary conditions. Lee \cite{Lee2012} studied mixed finite element methods with weak symmetry for the Zener, Kelvin-Voigt and Maxwell models, and employed the Crank-Nicolson scheme in the temporal discretization.  In \cite{Marques2012} Marques and Creuso gave an overview of numerical methods of viscoelasticity problems including finite element, boundary element and finite volume formulations. Li et al. \cite{Lihong2016} proposed a space-time continuous finite element method for a 2d viscoelastic wave equation.  In \cite{Wangshaojie2020}, Wang and Xie analyzed a  hybrid stress finite element method for the Maxwell model, where a second order implicit difference was used in the fully discrete scheme.  Recently,  Yuan and Xie \cite{Yuanhao}  showed that the mixed finite element framework  for Maxwell-model-based problems of wave propagation in linear viscoelastic solid  allows the use of a large class of existing mixed conforming finite elements for elasticity in the spatial discretization. 

This paper is to consider a class of weak Galerkin finite element discretizations of the quasistatic Maxwell viscoelastic model (\ref{modelproblem}).
 The weak Galerkin(WG)  method was firstly proposed  and analyzed by Wang and Ye  for second order elliptic problems \cite{Wangjunping2013A, Wangjunping2014}. Due to  adopting  weakly defined gradient/divergence operators over  functions  with discontinuity,  the WG method allows  in the finite element procedure the use of totally discontinuous functions and the local elimination of unknowns defined in the interior of elements.  Later on, this method was extended to some other models of  partial differential equations, such as
convection-diffusion equations \cite{Caowaixiang2021, Chengang2, Gaofuzheng2015,Lin2018, 
zhangtie2019}, linear elasticity problems \cite{Chengang1, Harper_Liu2019, Wangchunmei2016, Wang_Zhang2018}, Stokes equations \cite{Chen-Feng-Xie2016, Wang2014STOKES, zhai-Z-W2015, Zheng-Chen-Xie2017, Zheng-Xie2017}, Maxwell equations \cite{Mu-W-Y-Z2015, Shields2017}, natural convection problems \cite{Hanyihui2019B,Hanyihui2019A}, Biot models \cite{Chenyumei2016,Huxiaozhe2018}, 
 and biharmonic equations \cite{Burkardt-Gunz2020, Mulin2014A, Yexiu_Zhang2020}. 
 We  also refer the reader  to \cite{ Chen-W-W-Y2015, Libinjie2015A, Li-Xie2016, Li-Xie-Zhang2018} for  some fast solvers related to the WG algorithm.

In this contribution,  we develop semi-discrete and fully discrete WG methods for a velocity-stress system of the  quasistatic Maxwell viscoelastic model \eqref{modelproblem}, where the velocity variable $\bm{v}=\bm{u}_t$ is  introduced (cf. \eqref{trans_equation}). In the spatial discretization,  the  stress variable is approximated by  piecewise polynomials of degree  $k \ (k\geq 1)$,  the velocity variable is  approximated by  piecewise polynomials of degree   $k+1$, and the velocity trace   on the inter-element boundaries is  approximated by  piecewise polynomials of degree $k$. In the fully-discrete method, the backward Euler difference scheme is adopted for the temporal discretization. 

The rest of this paper is organized as follows. Section 2 introduces some notations and the weak variational problem. Section 3 and Section 4 are devoted to the  stability and error estimation  for the semi-discrete and fully discrete weak Galerkin schemes, respectively.   Finally, we report some numerical results to demonstrate the performance of the proposed WG  methods.

\section{  Weak  formulation}

We first introduce some notations. For any bounded domain $D\subset \mathbb{R}^s(s=d,d-1)$ and   nonnegative integer $m$, we denote by  $H^m(D)$ and $H_0^m(D)$  the usual $m$-th order Sobolev spaces   with norm $||\cdot||_{m,D}$ and semi-norm $|\cdot|_{m,D}$.   $H^0(D)=L^2(D)$ is  the space of square integrable functions defined on $D$. 
We use $(\cdot,\cdot)_{m,D}$ to denote the inner product of $H^m(D)$, with $(\cdot,\cdot)_D = (\cdot,\cdot)_{0,D}$.
When $D=\Omega$, we set $||\cdot||_m: = ||\cdot||_{m,\Omega}$, $|\cdot|_m:=|\cdot|_{m,\Omega}$ and $(\cdot, \cdot): = (\cdot,\cdot)_{\Omega}$.
In particular, for $D\subset \mathbb{R}^{d-1}$, we use $\langle\cdot, \cdot\rangle_{D}$ to replace $(\cdot, \cdot)_D$. 
For any integer $j\geq 0, P_j(D)$ denotes the set of all polynomials defined on $D$ with degree no greater than $j$.  

For  any vector-valued ( or tensor-valued)  space $X$, defined on $D$,  with norm $||\cdot||_X$,   we set
\begin{align*}
	L^p([0,T];X):=\left\{\boldsymbol{v}:[0,T]\rightarrow X;\ ||\boldsymbol{v}||_{L^p(X)}<\infty\right\},
\end{align*}
where
\begin{align*}
	||\boldsymbol{v}||_{L^p( X)}:=\left\{ \begin{array}{ll}
	(\int_{0}^{T}||\boldsymbol{v}(t)||_X^p)^{1/p} & \text{ if }1\leq p< \infty,\\
	\text{ess}\sup\limits_{0\leq t\leq T}||\boldsymbol{v}(t)||_X  & \text{ if } p=\infty,
	\end{array}
	\right.
\end{align*}	
and  $\boldsymbol{v}(t)$ abbreviates   $\boldsymbol{v}(\boldsymbol{x},t)$. 
For simplicity, we set $L^p(X):=L^p(0,T; X)$.  For any integer $r\geq0$, the spaces $H^r(X):=H^r(0,T; X)$ and   $C^r(X):=C^r([0,T];X)$  can be defined similarly.  

Let $\mathscr{T}_h=\bigcup \{ K \} $be a shape-regular   decomposition of the domain $\Omega\in \mathbb{R}^d(d=2,3)$  consisting of  
polygons/polyhedrons, in the sense that the following two assumptions hold (cf.\cite{Chengang1}):

\begin{quotation}
	(A1) There exists a positive constant $\theta_*$  such that  for each element $K\in \mathscr{T}_h$, there is a point $M_K\in K$ with $K$ being star-shaped with respect to every point in the ball of center $M_K$ and radius $\theta_* h_K$.
	\par 	(A2) There exists a positive constant $l_*$ such that for every element $K\in \mathscr{T}_h$, the distance between any two vertexes is no less than $l_*h_K$.
\end{quotation}
Let  $\mathcal{E}_h $ be the set of all edges/faces  of all elements in  $\mathscr{T}_h$.  
For any   $K\in \mathscr{T}_h$ and  $E \in \mathcal{E}_h$, we denote by $h_K$ and $h_E$ the diameters of $K$ and $E$, respectively, and set  $h:=\max\limits_{K\in\mathscr{T}_h}h_K$.  
Let $\nabla_h$ be the piecewise-defined gradient with respect to $\mathscr{T}_h$. 

For convenience,  throughout this paper we use $a \lesssim b$ to represent $a \leq Cb$, where  $C$ is a generic positive constant $C$ independent of the spatial mesh size $h$ and the temporal mesh size   $ \Delta t$.

Introducing the velocity variable $\bm{v}=\bm{u}_t$,  we  reformulate the quasistatic Maxwell viscoelastic model \eqref{modelproblem} as a velocity-stress form:

\begin{align}\label{trans_equation}
\left\{
\begin{array}{rl}
	-\mathbf{div}\bm\sigma =& \bm{f}(t) , \quad (x,t)\in \Omega \times[0,T], \\
	  \bm\sigma + \bm\sigma_t =& \mathbb{C}\bm{\varepsilon(v)}, \quad(x,t)\in \Omega\times[0,T],  \\
	 \bm{v}=&0,  \qquad  (x,t)\in \partial \Omega \times [0,T], \\
	\bm{\sigma}(0)=&\psi_0(x), \quad x\in \Omega.
\end{array}
\right.
\end{align}
It is easy to see that	 $\bm{u}(x,t)=\phi_0(x)+\int_0^t\bm{v}(x,s)ds.$
Define 
$$\bm L^2(\Omega,S):=\{\bm{\tau} = (\tau_{ij})_{d\times d} \in [L^2(\Omega)]^{d\times d}|\tau_{ij} = \tau_{ji} , \ i,j=1,2,\cdots, d\}.$$
Then, based on the system \eqref{trans_equation},   we can get the following  
weak problem:  
Find $(\bm{\sigma}, \bm{v})\in H^1(\bm L^2(\Omega,S))\times L^2([H_0^1(\Omega)]^d)$ such that

\begin{eqnarray}\label{original_formulation}
\left\{
\begin{array}{{rl}}
a(\bm{\sigma_t}, \bm{\tau}) + a(\bm\sigma,\bm\tau) - b(\bm\tau,\bm{v})=&0, \qquad  \forall   \bm{\tau}\in  \bm L^2(\Omega,S), \label{original_formulation1}\\
b( \bm\sigma, \bm{w})=& (\bm{f},\bm{w}),  \qquad  \forall \bm{w}\in  [H_0^1(\Omega)]^d,\label{original_formulation2}
\\
\bm\sigma(0) = &\psi_0(x), \qquad  x\in \Omega,
\end{array}
\right.
\end{eqnarray}
  for given
\begin{align}\label{cond-f}
\bm{f}\in H^1([L^2(\Omega)]^d),\quad  \psi_0\in  L^2(\Omega), 
\end{align}
where 
$a(\bm\sigma,\bm\tau):=(\mathbb{C}^{-1}\bm\sigma,\bm\tau)$ and 
$b(\bm\tau,\bm{w}):=(\bm{\tau},\bm{\varepsilon(\bm{w})}).$  

We need the following  continuous Gr\"onwall's inequality. 
\begin{Lem}Let $\phi(t)$ be such that
$$\dfrac{d\phi(t)}{dt}+\rho(t) \phi(t)\leq \psi(t), \hspace{2em}for \quad 0\leq t \leq T,$$
where $\rho(t), \psi(t) \in L^1([0,T])$. Then it holds 
\begin{align}
	\phi(t) \leq e^{-\int_0^t \rho(s)ds} \left(\phi(0)+\int_0^t\psi(s)e^{\int_0^s \rho(\tau)d\tau}ds\right), \hspace{2em} \forall t\in [0,T].\label{continuous_gronwall}
\end{align}
In particular, if $\rho \leq 0$ is a constant and $\psi(t)\geq 0$, then
	\begin{align}\label{Gr-2}
	\phi(t) \leq e^{-\rho T} \left(\phi(0)+\int_0^T\psi(s)ds\right), \hspace{2em} \forall t\in [0,T].
	\end{align}
\end{Lem}

	By following  a similar routine to that in    \cite{Rognes2010} for a weak formulation of the Maxwell model   with weak symmetry, we can derive    existence, uniqueness and stability results for the  system (\ref{original_formulation}):

\begin{Lem}
The weak problem \eqref{original_formulation} admits a unique solution $(\bm{\sigma},\bm{v}) \in H^1(\bm L^2(\Omega,S)) \times L^2([H^1_0(\Omega)]^d)$, and  the following stability results hold:
\begin{align}
||\bm\sigma(t)||_0^2&\lesssim  e^{-\frac{M_0}{M_1}t}||\psi_0||_0^2+\int_0^t e^{-\frac{M_0}{M_1}(t-s)}(||\bm{f}(s)||_0^2+||\bm{f}_t(s)||_0^2)ds,\label{stability1}\\
||\bm{v}(t)||_1^2 + ||\bm{\sigma_t}(t)||_0^2 &\lesssim ||\bm{\sigma}(t)||_0^2+||\bm{f}_t(t)||_0^2,\label{stability2}
\end{align}
for a.e. $t\in (0, T]$, where 
$M_0$ and $M_1$   are  positive constants given in  \eqref{assumed_condition}. 

\begin{proof}
On one hand, by the conditions \eqref{assumed_condition} and  \eqref{cond-f}
  there exists $\bm{\sigma}_e(t)\in  H^1(\bm L^2(\Omega,S)), \bm{v}_e(t)\in H^1([H^1_0(\Omega)]^d)$ solving the elasticity problem
\begin{align}
\left\{
\begin{array}{{rl}}
 a(\bm\sigma_e,\bm\tau) - b(\bm\tau,\bm{v}_e)=&0,  \forall   \bm{\tau}\in \bm L^2(\Omega,S), \\
b( \bm\sigma_e, \bm{w})=& (\bm{f},\bm{w}),  \forall \bm{w}\in [H_0^1(\Omega)]^d,
\end{array}
\right.
\end{align}
for $\forall t\in [0,T]$. 
Introduce 
$$\Sigma_0: = \{\bm{\tau}\in \bm L^2(\Omega,S)|b(\bm\tau,\bm{w}) = 0,\  \forall \bm{w}\in [H_0^1(\Omega)]^d\}.$$
From     \eqref{assumed_condition}  
we know that there exists  $\bm{\sigma}_0\in H^1(\Sigma_0)$ satisfying the ordinary differential equation
\begin{align}
\left\{ 
\begin{array}{rl}
a(\bm{\sigma}_{0,t}, \tau) + a(\bm{\sigma}_0, \bm{\tau}) =& -a(\bm{\sigma}_{e,t}, \tau), \qquad \forall \tau \in \Sigma_0, 
\\
\bm{\sigma_0}(0) = &\psi_0 - \bm{\sigma}_e(0).
\end{array}
\right.
\end{align} 

On the other hand, the Korn inequality indicates the inf-sup condition 
\begin{align}\label{continous_infsup}
|| \bm w||_1 \leq \beta ||\bm \varepsilon(\bm w)||_0\leq \beta\sup\limits_{0\neq \bm{\tau}\in \bm L^2(\Omega,S)} \dfrac{b(\bm{\tau},\bm{w})}{||\bm{\tau}||_0}, \ \forall \bm w\in [H_0^1(\Omega)]^d,
\end{align}
 with $\beta$ being a positive constant independent of  $\bm w $, which yields the existence of $\bm{v}_0(t)\in [H_0^1(\Omega)]^d$ for a.e. $t\in [0,T]$ such that
\begin{align}
a((\bm\sigma_{0}+\bm{\sigma}_{e})_t,\bm\tau) + a(\bm\sigma_{0},\bm\tau) - b (\bm{\tau}, \bm{v_0}) = 0, \quad  \forall   \bm{\tau}\in \bm L^2(\Omega,S).
\end{align}
As a result,  $\bm{\sigma} = \bm{\sigma}_0 + \bm{\sigma}_e$ and $ \bm{v} = \bm{v}_0 + \bm{v}_e $ solve the weak problem \eqref{original_formulation} for a.e. $t\in (0,T]$.

To prove the uniqueness of the solution, it suffices to establish the stability results \eqref{stability1} and \eqref{stability2}. 
To this end, we introduce an energy norm $||\cdot||_a$ on  $\bm L^2(\Omega,S)$ defined as $||\cdot||_a^2 := a(\cdot,\cdot)$. Then from  \eqref{assumed_condition} it follows 
\begin{align}\label{a-norm}
M_0||\bm\tau||_0^2\leq ||\cdot||_a^2\leq M_1||\bm\tau||_0^2, \quad \forall \bm\tau \in \bm L^2(\Omega,S).
\end{align}
We first prove \eqref{stability2}.  Take $\bm{\tau} = \bm{\sigma}_t$ in the first equation of \eqref{original_formulation} and differentiate the second equation of \eqref{original_formulation}  with respect to $t$ to obtain 
\begin{align}\label{y1}
a(\bm{\sigma},\bm{\sigma_t})+||\bm{\sigma_t}||_a^2 = (\bm{f}_t,\bm{v}).
\end{align}
In light of the first equation of \eqref{original_formulation}, the inf-sup condition \eqref{continous_infsup},  the Cauchy-Schwarz inequality and \eqref{a-norm}, 
we have 
\begin{align*} 
||\bm{v}(t)||_1 \leq \beta M_1^{1/2}(||\bm{\sigma}(t)||_a+||\bm{\sigma}_t(t)||_a),
\end{align*}
which,  together with \eqref{y1},  implies
\begin{align}\label{2.50}
||\bm{v}(t)||_1^2 + ||\bm{\sigma_t}(t)||_a^2 \leq  C||\bm{\sigma}(t)||_a^2 + ||\bm{f}_t(t)||_0^2.
\end{align}
Here $C$ is a positive constant depending only on $\beta, M_0,M_1$. Thus, from \eqref{a-norm} the desired estimate \eqref{stability2} follows.

The thing left is to show the stability \eqref{stability1}. 
Take $\bm{\tau} = \bm
{\sigma}$  and $\bm{w} = \bm{v}$ in  (\ref{original_formulation}) and employ the Young's inequality and \eqref{2.50} to get
\begin{align*}
	\dfrac{1}{2}\dfrac{d}{dt}||\bm{\sigma}(t)||_a^2+||\bm{\sigma}(t)||_a^2 &= (\bm{f},\bm{v})\leq \dfrac{q}{2}||\bm f(t)||_0^2+\dfrac{1}{2q}||\bm v(t)||_0^2 \\
	& \leq \dfrac{q}{2}||\bm f(t)||_0^2+\dfrac{C}{2q}(||\bm{\sigma}(t)||_a^2+||\bm{f}_t(t)||_0^2).
\end{align*}
Then, taking $q=\dfrac{M_1C}{2M_1-M_0}>0$ in this inequality implies 
\begin{eqnarray*}
	\dfrac{d}{dt}||\bm{\sigma}(t)||_a^2+\frac{M_0}{M_1}||\bm{\sigma}(t)||_a^2\leq c(||\bm{f}(t)||^2+||\bm{f}_t(t)||^2).
\end{eqnarray*}
Here $c$ is a positive constant depending only on $\beta, M_0,M_1$. 
Hence, using the Gr\"{o}nwall's inequality \eqref{continuous_gronwall}, we  obtain
\begin{align*}
	||\bm\sigma(t)||_a^2\leq e^{-\frac{M_0}{M_1}t}||\psi_0||_a^2+c\int_0^t e^{-\frac{M_0}{M_1}(t-s)}(||\bm{f}(s)||_0^2+||\bm{f}_t(s)||_0^2)ds,
\end{align*}
which, together with \eqref{a-norm}, yields \eqref{stability1}.
\end{proof}
\end{Lem}

\section{Semi-discrete Weak Galerkin  Method}
\subsection{Semi-discrete WG scheme}
We first follow \cite{Wangjunping2013A} to introduce the definitions of  discrete weak gradient/divergence operators.   

\begin{Def} 
For any $K\in \mathscr{T}_h$, $v\in \mathcal{V}(K):= \{  v=\{v_0,v_b\}:   v_0\in L^2(K),   v_b\in H^{\frac{1}{2}}(\partial K)     \}$ and   integer $j\geq 0$, the discrete weak gradient, $\nabla_{w,j,K} v\in [P_j(K)]^d$, of $v$   is defined by 
\begin{align} \label{DefinitionOfweakgradient}
( \nabla_{w,j,K} v,\bm{q} )_K:=-(v_0, \nabla \cdot \bm{q})_K + \langle v_b,\bm{q}\cdot \mathbf{n}_K \rangle_{\partial K}, \hspace{2mm} \forall \bm{q}\in [P_j(K)]^d,
\end{align}
where $\mathbf{n}_K$ is the unit outward normal vector along $\partial K$. The global discrete weak gradient operator $\nabla_{w,j}$ on $\mathcal{V}(\mathscr{T}_h):=\{v:\ v|_K\in \mathcal{V}(K), \ \forall K\in \mathscr{T}_h\}$ is defined  by
$$\nabla_{w,j}|_K=\nabla_{w,j,K}, \forall K\in \mathscr{T}_h.$$ 
 For a vector $\bm v = (v_1, \cdots , v_d)^T \in [\mathcal{V}(\mathscr{T}_h)]^d$,
  its discrete weak gradient $\nabla_{w,j}\bm v$  is defined as
$$ \nabla_{w,j}\bm v:=(\nabla_{w,j} v_1, \cdots, \nabla_{w,j} v_d)^T.$$
\end{Def}

\begin{Def}\label{weak-div} For any $K\in \mathscr{T}_h$, $\bm{v}\in \mathcal{W}(K):=\{\bm{v}=\{\bm{v_0},\bm{v_b}\}:\bm{v_0}\in [L^2(K)]^d, \bm{v_b}\cdot\bm{n}_K\in H^{-1/2}(\partial K)\}$ and   integer $j\geq 0$,
the discrete weak divergence, $\nabla_{w,j,K}\cdot \bm{v}\in P_j(K)$, of $\bm v$ is defined by
\begin{align}
(\nabla_{w,j,K}\cdot \bm{v}, q)_K = -(\bm{v_0},\nabla q)_K + \langle \bm{v_b}\cdot \bm{n}_K, q\rangle_{\partial K}, \quad \forall q\in P_j(K).
\end{align}
The  global discrete weak divergence operator $\nabla_{w,j}\cdot$ is defined by $$\nabla_{w,j}\cdot|_K = \nabla_{w,j,K}\cdot, \quad \forall  K\in \mathscr{T}_h. $$
\end{Def}

For any  $K\in \mathscr{T}_h$, $E\in \mathcal{E}_h$ and any   integer $j\geq 0$,  let 
$$Q_j^0: L^2(K)\rightarrow P_{j}(K), \quad 
Q_j^b: L^2(E)\rightarrow P_j(E)$$
be the usual  $L^2$ projection operators.  For convenience, vector and tensor analogues of ${Q}_j^0$ and ${Q}_j^b$ are still denoted by $ {Q}_j^0$ and $ {Q}_j^b$, respectively.  
 
 For any integer $k\geq 1$,  we introduce the following finite dimensional spaces:
 \begin{align}
\Sigma_{h}:&=\{\bm{\tau_{h}}\in \bm L^2(\Omega,S):   \bm{\tau_{h}}|_K\in [P_k(K)]^{d\times d}, \hspace{1mm} \forall  K\in \mathscr{T}_h\} ,
 \\
 V_{h}:&=\{\bm{v_h}=\{\bm{v_{h0},v_{hb}}\} :\bm{v_{h0}}|_K\in [P_{k+1}(K)]^d, \bm{v_{hb}}|_E\in [P_k(E)]^d,\hspace{1mm}\forall  K\in \mathscr{T}_h, E\in \mathcal{E}_h\},\\
  V_h^0:&= \{\bm{v_h}\in V_h:\bm{v_{hb}}|_{\partial \Omega}=0\}.
 \end{align}

The semi-discrete WG scheme reads as follows: For any $t\in [0,T]$, find $\bm{\sigma_{h}}(\cdot, t)\in \Sigma_h$, $\bm{v_h}(\cdot,t) = \{\bm{v_{h0}}(\cdot,t), \bm{v_{hb}}(\cdot,t)\}\in V_h^0$ such that
\begin{subequations}\label{semidiscretescheme}
\begin{align}
a_h(\bm{\sigma_{h,t},\tau_h}) + a_h(\bm{\sigma_{h},\tau_h})  -b_h(\bm{\tau_h,v_h})=0 ,& \forall \bm{\tau_h}\in \Sigma_h \label{semidiscretescheme_1},
\\
b_h( \bm{\sigma_{h},w_h}) + s_h(\bm{v_h,w_h}) =(\bm{f, w_{h0}}), & \forall  \bm{w_h}\in V_h^0\label{semidiscretescheme_2}
\\
 \bm{\sigma_{h}}(0)= {Q}_k^0\psi_{0}, & \label{semidiscretescheme_3}
\end{align}
\end{subequations}
where
\begin{align*}
	&a_h(\bm{\sigma_h,\tau_h})=(\mathbb{C}^{-1}\bm{\sigma_h,\tau_h}), \qquad 
	  b_h(\bm{\tau_h,  w_h}) = (\bm{\varepsilon_{w,k}( w_h), \tau_h}) , 
	\\
	&s_h(\bm{v_h,w_h}) = \langle \alpha(Q_k^b\bm{v_{h0}-v_{hb}}), Q_k^b\bm{w_{h0}-w_{hb}}\rangle_{\partial \mathscr{T}_h},
	\end{align*}
with  $ \bm{\varepsilon_{w,k}(w_h)}: = \left({\bm{\nabla_{w,k} w_h+(\nabla_{w,k} w_h)^T}}\right)/{2}$, $ \langle \cdot, \cdot\rangle_{\partial \mathscr{T}_h}:=\sum_{K\in\mathscr{T}_h} \langle \cdot, \cdot\rangle_{\partial K}$, and 
the stabilization parameter $\alpha|_E =   h_E^{-1}$ for any $E\in \mathcal{E}_h$. 

\begin{Rem}Notice that by the definition of the discrete weak gradient  we have 
\begin{equation}\label{bh-new}
b_h(\bm{\tau_h,  w_h}) =  (\bm{\nabla_{w,k} w_h, \tau_h})=- (\bm{ w_{h0}, \nabla\cdot\tau_h})_{\mathscr{T}_h}+
\langle \bm {w_{hb}, \tau_h n}\rangle_{\partial \mathscr{T}_h}.
\end{equation}
Then the equations  \eqref{semidiscretescheme_1}-\eqref{semidiscretescheme_2} lead to  the relations
\begin{subequations}\label{hdgscheme}
	\begin{align}
	&a_h(\bm{\sigma_{h,t},\tau_h}) + a_h(\bm{\sigma_{h},\tau_h})  + (\bm{v}_{h0},\nabla\cdot \bm\tau_h) - \langle \bm{v}_{hb},\bm\tau_h \bm n\rangle_{\partial \mathscr{T}_h}=0, \label{hdgscheme1} \\
	&  - (\bm{\nabla\cdot \sigma_{h}}, \bm w_{h0}) + \langle \alpha(Q_k^b\bm{v_{h0}-v_{hb}}),\bm{w_{h0}}\rangle_{\partial \mathscr{T}_h} = (\bm{f,\bm w_{h0}}), \label{hdgscheme2}
	\\
	& \langle \bm{\sigma_{h} n} - \alpha(Q_k^b \bm{v_{h0}-v_{hb}}),\bm{w_{hb}} \rangle_{\partial \mathscr{T}_h}=0, \label{hdgscheme3}
	\end{align}
\end{subequations}
  for all $(\bm{\tau_h,\{w_{h0},w_{hb}}\})\in \Sigma_h\times V_h^0$. 

\end{Rem}

By using standard techniques, we can show the   existence and uniqueness of the semi-discrete solution.
\begin{thm} \label{exist_unique_semi}
The semi-discrete scheme (\ref{semidiscretescheme}) admits  a unique solution $(\bm{\sigma_{h},  v_h})\in \Sigma_h\times V_h^0.$
\begin{proof}
Let $ \{\Phi_i\}_{i=1}^{r_1}$ and $ \{\{\phi_{0i}\}_{i=1}^{r_2},  \{\phi_{bi}\}_{i=1}^{r_3}\}$ be the  basis functions of $\Sigma_{h}$ and $V_{h}^0$, respectively. 
We write $$\bm{\sigma_{h}}(t)=\sum_{i=1}^{r_1}\eta_i(t)\Phi_i ,\quad \bm{v_{h0}}(t) = \sum_{i=1}^{r_2}\beta_i(t)\phi_{0i},\quad \bm{v_{hb}} =\sum_{i=1}^{r_3} \gamma_i(t)\phi_{bi},\quad \mathcal{F}_i= (f,\phi_{0j}) ,$$
and denote by $\eta(t),\beta(t),\gamma(t)$ the corresponding vectors of $\eta_i(t),\beta_i(t),\gamma_i(t)$, respectively. Let $\mathcal{M}_{s,ij}$ the $(i,j)-th$ components of  matrix $\mathcal{M}_s(s=0,1,\cdots,6)$ be given by
\begin{align*}
& \mathcal{M}_{0,ij}=( \mathbb{C}^{-1}\Phi_j, \Phi_i ),  
&& \mathcal{M}_{1,ij} =( \phi_{0j},\nabla\cdot\Phi_i)      
			\\
& \mathcal{M}_{2,ij}=  -\langle \phi_{bj},  \Phi_i \bm{n}\rangle_{\partial \mathscr{T}_h}, 
&&\mathcal{M}_{3,ij}= \langle\alpha Q_k^b\phi_{0j},\phi_{0i}\rangle_{\partial \mathscr{T}_h}, 
\\
&\mathcal{M}_{4,ij}= -\langle\alpha \phi_{bj},  \phi_{0i}\rangle_{\partial \mathscr{T}_h}    ,
&&\mathcal{M}_{5,ij}=-\langle\alpha Q_k^b\phi_{0j}, \phi_{bi})\rangle_{\partial \mathscr{T}_h}
			\\
&\mathcal{M}_{6,ij}=\langle\alpha \phi_{bj},\phi_{bi}\rangle_{\partial \mathscr{T}_h}. 
& &
\end{align*}
Then  the system (\ref{semidiscretescheme}) can be written as the following matrix  forms: 
\begin{align}
\mathcal{M}_0\dfrac{d{\eta}(t)}{dt} + \mathcal{M}_0 {\eta}(t)+ \mathcal{M}_1{\beta}(t) + \mathcal{M}_2{\gamma}(t)& = 0,\label{q1}
 \\
-\mathcal{M}_1^T{\eta}(t)+\mathcal{M}_3{\beta}(t)
+ \mathcal{M}_4{\gamma}(t)& = \mathcal{F}(t),\label{q2}
 \\
-\mathcal{M}_2^T{\eta}(t) + \mathcal{M}_5\beta(t) + \mathcal{M}_6\gamma(t) &= 0.\label{q3}
\end{align}
Here we have used the relation \eqref{bh-new} for the terms $b_h(\cdot,\cdot)$ in the scheme. 
 Since $\mathcal{M}_0, \mathcal{M}_6$ are symmetric positive defined, we can eliminate $\beta(t)$ and $\gamma(t)$ from (\ref{q1})-(\ref{q3}) to get
\begin{align}
\mathcal{M}_0\dfrac{d\eta(t)}{dt} + \mathcal{P}\eta(t) = \mathcal{Q}(t),
\end{align}
where
\begin{align*}
\mathcal{P} :=& \mathcal{M}_0+\mathcal{M}_2\mathcal{M}_6^{-1}\mathcal{M}_2^T+(\mathcal{M}_1-\mathcal{M}_2\mathcal{M}_6^{-1}\mathcal{M}_5)(\mathcal{M}_3-\mathcal{M}_4\mathcal{M}_6^{-1}\mathcal{M}_5)^{-1}(\mathcal{M}_1^T-\mathcal{M}_4\mathcal{M}_6^{-1}\mathcal{M}_2^T),
\\
\mathcal{Q} :=& (\mathcal{M}_2\mathcal{M}_6^{-1}\mathcal{M}_5-\mathcal{M}_1)(\mathcal{M}_3-\mathcal{M}_4\mathcal{M}_6^{-1}\mathcal{M}_5)^{-1}\mathcal{F}(t).
\end{align*}	
By the standard  theory of ordinary differential equations (cf. \cite{Earl1955}), the above system (also the scheme (\ref{semidiscretescheme})),  has a unique solution. This completes the proof.
\end{proof}	
\end{thm}

\subsection{A priori error estimation}

To establish error estimates for the proposed WG scheme, we need the following properties of 
  the $L_2$
-projections $Q_{j}^0 ,Q_j^b$ with nonnegative integer j.

\begin{Lem}\cite{Chengang1}\label{Communityproperty} It holds the commutative property
\begin{align}
\nabla_{w,j}\{Q_{j+1}^0\bm v,Q_j^b\bm v\}   =  {Q}_j^0\nabla \bm v , for \hspace{2mm}all\hspace{2mm}\bm{v}\in[H^1(K)]^d ,
\end{align}
\end{Lem}

\begin{Lem}\cite{Chengang1,Shizhongci2013}\label{projection_property} Let $m$ be an integer with $1\leq m \leq j+1$. For any $K\in \mathscr{T}_h, E\in\mathcal{E}_h$, it holds
\begin{subequations}
\begin{align}
& ||v-Q_j^0v||_{0,K}+h_K|v-Q_j^0v|_{1,K} \lesssim h_K^m |v|_{m,K}, && \forall v\in H^m(K), \\
& ||v-Q_j^bv||_{0,\partial K} \lesssim h_k^{m-1/2}|v|_{m,K},   && \forall v\in H^m(K), \\
& |v-Q_j^0v|_{s,K} \lesssim h_K^{m-s}|v|_{m,K}, && \forall v\in H^m(K), 0\leq s\leq m, \\
&||Q_j^0 v||_{0,K}\leq ||v||_{0,K} , && \forall v\in L^2(K), \\
& ||Q_j^b v||_{0,E} \leq ||v||_{0,E}, && \forall v\in L^2(E) , 
\end{align}
\end{subequations}
\end{Lem}

For the the bilinear forms $a_h(\cdot,\cdot)$ and  $b_h(\cdot,\cdot)$, we   easily get the following continuity  and coercivity results.
\begin{Lem} For all $\bm{\sigma_{h}, \tau_{h}}\in \Sigma_h,  \bm{v_h}=\{\bm{v_{h0},v_{hb}}\}\in V_h$,  it holds 
\begin{align}
a_h(\bm{\sigma_{h0},\tau_{h0}})&\leq M_1||\bm{\sigma_{h0}}||_0||\bm{\tau_{h0}}||_0,\\
b_h(\bm{\tau_{h}},\bm{w_h})&\leq  ||\bm{\tau_{h}}||_0||\bm{\varepsilon_{w,k}(w_h)}||_0,\\
 a_h(\bm{\tau_h,\tau_h})&\geq M_0\|\bm{\tau_h}\|_0^2.
\end{align}
\end{Lem}	
We also need the following \emph{inf-sup} stability condition for the bilinear form $b_h(\cdot,\cdot)$.
\begin{Lem} \label{infsupcondition}
	{\rm\cite{Chengang1}}For any $\bm{w_h}=\{\bm{w_{h0},w_{hb}}\}\in V_h$, it holds
	\begin{align}\label{inf-sup-dis}
	 ||\bm{\varepsilon_{w,k}(w_h)}||_0\lesssim\sup_{0\neq \bm{\tau_{h}}\in \Sigma_{h}}\dfrac{b_h(\bm{\tau_{h}},\bm{w_h})}{||\bm{\tau_{h}}||_0} .
	\end{align}
\end{Lem}

\begin{Lem}\label{simp}{\rm\cite{Chengang1}}
For any $\bm{w_h }= \{\bm{w_{h0}}, \bm{w_{hb}}\}\in V_h^0$ and sufficiently small $h$, it holds
\begin{align}
||\bm{\nabla_h w_{h0}}||_0^2 &\lesssim ||\bm{\varepsilon_h(w_{h0})}||_0^2 + ||\alpha^{1/2}\bm{(Q_k^b w_{h0} - w_{hb})}||_{\partial \mathscr{T}_h}^2 , \label{s1}
\\
||\bm{\varepsilon_h(w_{h0})}||_0^2&\lesssim ||\bm{\varepsilon_{w,k}(w_h)}||_0^2 + ||\alpha^{1/2}\bm{(Q_k^bw_{h0}-w_{hb})}||_{\partial \mathscr{T}_h},\label{s2}
\end{align}
where $\bm{\varepsilon_h(w)}:=\left(\bm{\nabla_h w + (\nabla_h w)^T}\right)/2$ and $||\cdot||_{\partial \mathscr{T}_h}:=\langle \cdot, \cdot\rangle_{\partial \mathscr{T}_h}^{1/2}.$
\end{Lem}


The following lemma shows  the error equations of the weak solution $(\bm{\sigma}, \bm{v})$ and its projection $( {Q}_k^0\bm{\sigma}, \{ {Q}_{k+1}^0\bm v,  {Q}_k^b\bm v\})$.
\begin{Lem} \label{semi_discrete_proj_problem}
Let $(\bm{\sigma}, \bm{v})\in H^1(\bm L^2(\Omega,S)\bigcap  \bm H ({\bf div},\Omega))\times L^2([H_0^1(\Omega)]^d)$
be the weak solution of system \eqref{original_formulation}, then, for all $\bm{\tau_h}\in\Sigma_{h}$ and $\bm{w_h=\{w_{h0},w_{hb}\}}\in V_h^0$ it holds
\begin{subequations} \label{semi_discrete_proj_scheme}
\begin{align}
&a_h({Q}_k^0\bm{\sigma}_t ,\bm{\tau}_h) + a_h( {Q}_k^0 \bm{\sigma},\bm{\tau}_h) - b_h(\bm{\tau}_h,\{{Q}_{k+1}^0\bm v, {Q}_k^b\bm v\})=a_h(({Q}_k^0\bm{\sigma - \sigma})_t,\bm{\tau}_h), \label{semi_discrete_proj_scheme_1}
\\
&b_h ( {Q}_k^0\bm{\sigma},\bm w_h) + s_h(\{{Q}_{k+1}^0\bm v, {Q}_k^b\bm v\}, \bm w_h) = (\bm{f},\bm w_{h0}) + l_1(\bm{\sigma},\bm w_h)+l_2(\bm{v},\bm w_h),\label{semi_discrete_proj_scheme_2}
\end{align}
\end{subequations}
where 
$$l_1(\bm{\sigma},\bm w_h) := \langle \bm w_{h0}-\bm w_{hb},\bm{\sigma n} - {Q}_k^0 \bm{\sigma n}\rangle_{\partial \mathscr{T}_h},\quad 
l_2(\bm{v},\bm w_h):=\langle \alpha(Q_k^bQ_{k+1}^0\bm{v}-Q_k^b\bm{v}), Q_k^b\bm w_{h0} - \bm w_{hb}\rangle_{\partial \mathscr{T}_h}.$$
 
\begin{proof}
By  the commutative property (\ref{Communityproperty}) and the definitions of $a_h(\cdot,\cdot)$ and $ b_h(\cdot,\cdot)$, we   obtain
\begin{align}
&a_h({Q}_k^0\bm{\sigma}_t,\tau_h)+a_h({Q}_k^0\bm{\sigma}, \tau_h) - b_h(\tau_h, \{{Q}_{k+1}^0\bm v, {Q}_k^b\bm v\})
\nonumber \\
=& (\mathbb{C}^{-1}{Q}_k^0\bm{\sigma}_t, \tau_h) + (\mathbb{C}^{-1}{Q}_k^0\bm{\sigma}, \tau_h) - (\bm{\nabla_{w,k}} \bm{Q}_h^k\bm{v},\tau_h) 
\nonumber \\
 = &
 (\mathbb{C}^{-1}\bm{\sigma}_t, \tau_h) + (\mathbb{C}^{-1}{Q}_k^0\bm{\sigma}, \tau_h) - ({Q}_k^0\nabla \bm{v},\tau_h)+(\mathbb{C}^{-1}({Q}_k^0\bm{\sigma-\sigma})_t, \tau_h)
\nonumber \\ 
 = & (\mathbb{C}^{-1}\bm{\sigma}_t, \tau_h) + (\mathbb{C}^{-1}\bm{\sigma}, \tau_h) - (\nabla \bm{v},\tau_h)+(\mathbb{C}^{-1}({Q}_k^0\bm{\sigma-\sigma})_t, \tau_h)
\nonumber \\
=&(\mathbb{C}^{-1}({Q}_k^0\bm{\sigma-\sigma})_t, \tau_h) .
\end{align}
From the definition of weak gradient, the projection property  and the Green's formula, it follows
\begin{align}
&b_h ( {Q}_k^0\bm{\sigma},w_h) + s_h(\{{Q}_{k+1}^0\bm v, {Q}_k^b\bm v\}, w_h) 
\nonumber \\
=&	
 (\bm{\nabla_{w,k}} w_h, {Q}_k^0\bm{\sigma})+\langle \alpha(Q_k^bQ_{k+1}^0v - Q_k^bv), Q_k^b w_{h0} - w_{hb} \rangle_{\partial \mathscr{T}_h}  
 \nonumber \\
 = & 
 - (w_{h0}, \nabla_h\cdot {Q}_k^0\bm{\sigma}) + \langle w_{hb}, {Q}_k^0 \bm{\sigma} n\rangle_{\partial \mathscr{T}_h} + \langle \alpha(Q_k^bQ_{k+1}^0\bm{v} - Q_k^b\bm{v}), Q_k^b w_{h0} - w_{hb} \rangle_{\partial \mathscr{T}_h}  \nonumber \\
 = & 
(\nabla_h w_{h0},  {Q}_k^0\bm{\sigma}) - \langle w_{h0}-w_{hb}, {Q}_k^0 \bm{\sigma} n\rangle_{\partial \mathscr{T}_h} + \langle \alpha(Q_k^bQ_{k+1}^0\bm{v} - Q_k^b\bm{v}), Q_k^b w_{h0} - w_{hb} \rangle_{\partial \mathscr{T}_h}  \nonumber \\
 = &
(-\nabla\cdot\bm{\sigma}, w_h) + \langle w_{h0}-w_{hb}, (\bm{\sigma} -{Q}_k^0 \bm{\sigma} )n\rangle_{\partial \mathscr{T}_h} + \langle \alpha(Q_k^bQ_{k+1}^0\bm{v} - Q_k^b\bm{v}), Q_k^b w_{h0} - w_{hb} \rangle_{\partial \mathscr{T}_h} \nonumber \\
 = & (\bm{f},w_h) + l_1(\bm{\sigma}, w_h) + l_2(\bm{v},w_h).
\end{align}
This finishes the proof.
\end{proof}
\end{Lem}

\begin{Lem}\label{approximationlemma}
Let $ (\bm{\sigma}, \bm{v})\in  H^1( \bm L^2(\Omega,S)\bigcap[H^{k+1}(\Omega)]^{d\times d})\times L^2([H_0^1(\Omega)\bigcap H^{k+2}(\Omega)]^d )$ be the weak solution of system \eqref{original_formulation} and $w_h=\{w_{h0},w_{hb}\}\in V_h$, it holds
\begin{align}
|l_1(\bm{\sigma}, w_h)| \lesssim & h^{k+1}|\bm{\sigma}|_{k+1}||\nabla_h w_{h0}||_0 + h^{k+1}|\bm{\sigma}|_{k+1}||\alpha^{1/2}(Q_k^bw_{h0}-w_{hb})||_{\partial  \mathscr{T}_h},  
	\\
|l_2(\bm{v}, w_h)| \lesssim& h^{k+1}|\bm{v}|_{k+2}||\alpha^{1/2}(Q_k^bw_{h0} - w_{hb})||_{\partial  \mathscr{T}_h} .
\end{align}
	
\begin{proof}
%
Using the Cauchy-Schwarz inequality, the projection properties, the trace inequality and the  triangle inequality, we  obtain
\begin{align}
|l_1(\bm{\sigma},w_h)|\leq &||w_{h0}-w_{hb}||_{\partial  \mathscr{T}_h}||\bm{\sigma} \bm{n} - Q_k^0 \bm{\sigma} \bm{n}||_{\partial \mathscr{T}_h} 
\nonumber \\
 = & ||\alpha^{1/2}(w_{h0}-w_{hb})||_{\partial  \mathscr{T}_h}||\alpha^{-1/2}(\bm{\sigma} \bm{n} - Q_k^0 \bm{\sigma} \bm{n})||_{\partial \mathscr{T}_h} 
\nonumber \\
 \lesssim & h^{k+1}|\bm{\sigma}|_{k+1} ||\alpha^{1/2}(w_{h0}-w_{hb})||_{\partial  \mathscr{T}_h} 
\nonumber \\
 \leq & 
h^{k+1}|\bm{\sigma}|_{k+1} ||\alpha^{1/2}(w_{h0}-Q_k^bw_{h0})||_{\partial  \mathscr{T}_h} + h^{k+1}|\bm{\sigma}|_{k+1} ||\alpha^{1/2}(Q_k^bw_{h0}-w_{hb})||_{\partial  \mathscr{T}_h}
\nonumber \\
 \lesssim & h^{k+1}|\bm{\sigma}|_{k+1}h^{-1/2}||w_{h0}-Q_k^bw_{h0}||_{\partial \mathscr{T}_h} + h^{k+1}|\bm{\sigma}|_{k+1}||\alpha^{1/2}(Q_k^bw_{h0}-w_{hb})||_{\partial  \mathscr{T}_h}
\nonumber \\
 \lesssim & 
h^{k+1}|\bm{\sigma}|_{k+1}||\nabla_h w_{h0}||_0 + h^{k+1} |\bm{\sigma}|_{k+1}||\alpha^{1/2}(Q_k^bw_{h0}-w_{hb})||_{\partial  \mathscr{T}_h} .
\end{align}
		
Similarly, by the Cauchy-Schwarz inequality and the projection properties  we  get
\begin{align}
|l_2(\bm{v},w_h)|\leq & ||\alpha^{1/2}(Q_k^bQ_{k+1}^0\bm{v}-Q_k^b\bm{v})||_{\partial  \mathscr{T}_h} ||\alpha^{1/2}(Q_k^bw_{h0} - w_{hb})||_{\partial  \mathscr{T}_h}
\nonumber \\
 \lesssim& ||\alpha^{1/2}(Q_{k+1}^0\bm{v-v})||_{\partial  \mathscr{T}_h} ||\alpha^{1/2}(Q_k^bw_{h0} - w_{hb})||_{\partial  \mathscr{T}_h}
\nonumber \\
 \lesssim & h^{k+1}|\bm{v}|_{k+2}||\alpha^{1/2}(Q_k^bw_{h0} - w_{hb})||_{\partial  \mathscr{T}_h}.
\end{align}
This completes the proof.
\end{proof}
\end{Lem}



The following lemma gives  an   estimate of the error  
between  the semi-discrete solution $(\bm{\sigma}_{h}, \bm v_h=\{\bm{v}_{h0},\bm{v}_{hb}\}) $ and the projection $( {Q}_k^0\bm{\sigma}, \{ {Q}_{k+1}^0\bm v, {Q}_k^b\bm v\})$ of the weak solution.

\begin{Lem}\label{approximation_error} Let $ (\bm{\sigma}, \bm{v})\in  C^1( \bm L^2(\Omega,S)\bigcap[H^{k+1}(\Omega)]^{d\times d})\times C^1([H_0^1(\Omega)\bigcap H^{k+2}(\Omega)]^d )$ be the weak solution of   \eqref{original_formulation} and $(\bm{\sigma}_{h}, \{\bm{v}_{h0},\bm{v}_{hb}\})\in C^1(\Sigma_{h})\times C^1(V_h^0)$ be the  semi-discrete solution of the WG scheme \eqref{semidiscretescheme}.  
Then   it holds
\begin{align}
||\zeta_{h}||_0^2+s_h(\xi_h, \xi_h )
 \lesssim& h^{2k+2}\left( \tilde M_0(\sigma, v)+\tilde M_2(\sigma, v)\right), \label{error_stress_est}\\
  ||\bm{\varepsilon_h}(\xi_{h0})||_0^2 \lesssim&   h^{2k+2}\left(\tilde M_0(\sigma, v)+\tilde M_1(\sigma, v)+\tilde M_2(\sigma, v)\right), \label{error_strain_est}
\end{align}
where $\zeta_{h} := {Q}_k^0 \bm{\sigma} - \bm{\sigma_{h}}$, $\xi_h:=\{\xi_{h0}, \xi_{hb}\}$ with 
$\xi_{h0} = Q_{k+1}^0\bm{v-v}_{h0}$ and $\xi_{hb} = Q_k^b\bm{v-v}_{hb}$, 
\begin{align*}
\tilde M_0(\sigma, v):=& |\bm{\sigma}(0)|_{k+1}^2 + |\bm v(0)|_{k+2}^2 + |\bm{\sigma}_t(0)|_{k+1}^2,\\
 \tilde M_1(\sigma, v):=&|\bm{\sigma}_t|_{k+1}^2+|\bm{v}_t|_{k+2}^2,\\
\tilde M_2(\sigma, v):=&\int_0^t(|\bm{\sigma}|_{k+1}^2+|\bm{v}|_{k+2}^2+|\bm{\sigma}_t|_{k+1}^2+|\bm{v}_t|_{k+1}^2) ds.
\end{align*}
\begin{proof} 

Substract (\ref{semidiscretescheme_1}) and (\ref{semidiscretescheme_2}) from  (\ref{semi_discrete_proj_scheme_1}) and (\ref{semi_discrete_proj_scheme_2}), respectively, we obtain
\begin{align}
&a_h(\zeta_{h,t},\tau_h)+a_h(\zeta_{h},\tau_h)-b_h(\tau_h,\xi_h) = a_h((Q_k^0\bm{\sigma} - \bm{\sigma})_t,\tau_h),  \label{s6}
		\\
&b_h(\zeta_{h}, w_h) + s_h(\xi_h, w_h) = l_1(\bm{\sigma},w_h) + l_2(\bm{v},w_h). \label{s7}
\end{align}
Taking $(\tau_h,w_h)=(\tau_h,\{w_{h0},w_{hb}\})=(\zeta_{h}, \{\xi_{h0} , \xi_{hb} \}) = (\zeta_{h}, \xi_h )$ in the above equations yields 
\begin{align}
\dfrac{1}{2}\dfrac{d}{dt}a_h(\zeta_{h},\zeta_{h}) + a_h(\zeta_{h},\zeta_{h}) +s_h(\xi_h, \xi_h )= a_h(({Q}_k^0\bm{\sigma - \sigma})_t,\zeta_{h})+  l_1(\bm{\sigma},\xi_h ) + l_2(\bm{v},\xi_h ).
\end{align}
From Lemma \ref{approximationlemma}, Lemma \ref {simp} and the Young's inequality with any $\kappa >1$ it follows
\begin{align*}
&\dfrac{1}{2}\dfrac{d}{dt}||\zeta_{h}||_a^2+||\zeta_{h}||_a^2 +s_h(\xi_h, \xi_h)
\nonumber \\ 
 \leq & \dfrac{1}{2}||({Q}_k^0\bm{\sigma - \sigma})_t||_a^2+\dfrac{1}{2} ||\zeta_h||_a^2+ Ch^{k+1}|\bm{\sigma}|_{k+1}||\nabla_h w_{h0}||_0+Ch^{k+1}|\bm{\sigma}|_{k+1}||\alpha^{1/2}(Q_k^b\xi_{h0}-\xi_{hb})||_{\partial  \mathscr{T}_h}
\nonumber \\
&+ h^{k+1}|\bm{v}|_{k+2}||\alpha^{1/2}(Q_k^b\xi_{h0} - \xi_{hb})||_{\partial  \mathscr{T}_h}
\nonumber \\
 \leq & \dfrac{1}{2M_0}||({Q}_k^0\bm{\sigma - \sigma})_t||_0^2+\dfrac{1}{2} ||\zeta_h||_a^2+ \kappa Ch^{2k+2}(|\bm{\sigma}|_{k+1}^2+|\bm{v}|_{k+2}^2) +  \dfrac{C}{2\kappa}||\nabla_h \xi_{h0}||_0^2+\dfrac{C}{\kappa}||\alpha^{1/2}(Q_k^b\xi_{h0}-\xi_{hb})||_{\partial  \mathscr{T}_h}^2
\nonumber \\
\leq & \dfrac{C}{2M_0}h^{2k+2}|\bm{\sigma}_t|_{k+1}^2+\dfrac{1}{2}||\zeta_{h}||_a^2+  Ch^{2k+2}(|\bm{\sigma}|_{k+1}^2+ |\bm v|_{k+2}^2)+ \dfrac{C}{\kappa}||\alpha^{1/2}(Q_k^b\xi_{h0}-\xi_{hb})||_{\partial \mathscr{T}_h}^2  + \dfrac{C}{2\kappa}||\bm{\varepsilon_{h}} (\xi_h)||_0^2,  
\end{align*}
which implies that
\begin{align}\label{x1}
\dfrac{d}{dt}||\zeta_{h}||_a^2+||\zeta_{h}||_a^2 +s_h(\xi_h, \xi_h)
\lesssim  h^{2k+2} (|\bm{\sigma}|_{k+1}^2+ |\bm v|_{k+2}^2+|\bm{\sigma}_t|_{k+1}^2)  + \dfrac{C}{\kappa}||\bm{\varepsilon_{w,k}}(\xi_h)||_0^2.
\end{align}

By Lemma \ref{infsupcondition}, Lemma \ref{projection_property} and equation \eqref{s6},  we have
\begin{align}\label{strain_estim}
||\bm{\varepsilon_{w,k}}(\xi_h)||_0 \lesssim & \sup_{0\neq \tau_{h}\in \Sigma_{h}}\dfrac{b_h(\tau_h, \xi_h)}{||\tau_{h}||_0} 
		\nonumber \\
=&\sup_{0\neq \tau_{h}\in \Sigma_{h}}\dfrac{a_h(\zeta_{h,t}, \tau_h) + a_h(\zeta_{h}, \tau_h)-a_h((Q_k^0\bm{\sigma}-\bm{\sigma})_t, \tau_h) }{||\tau_{h}||_0}
		\nonumber \\
\leq &  c( ||\zeta_{h,t}||_0 + ||\zeta_{h}||_0+h^{k+1}|\bm{\sigma}_t|_{k+1}).
\end{align}
Here $c$ is a positive constant independent of $h$. 
To bound the term $||\zeta_{h,t}||_0$, substitute $\tau_h = \zeta_{h,t}$ into   (\ref{s6}) and take $w_h = \xi_h$ in (\ref{s7}) after differentiating in time, then we get
\begin{align*}
a_h(\zeta_{h,t},\zeta_{h,t})+a_h(\zeta_{h},\zeta_{h,t})-b_h(\zeta_{h,t},\xi_h)=&a_h(({Q}_k^0\bm{\sigma - \sigma})_t, \zeta_{h,t}),
\\
b_h(\zeta_{h,t}, \xi_h) + \dfrac{1}{2}\dfrac{d}{dt}s_h(\xi_h,\xi_h) =& l_1(\bm{\sigma}_t, \xi_h) + l_2(\bm{v}_t, \xi_h).
\end{align*}
Summing up the above two equalities and using Lemma \ref{approximationlemma}, Lemma \ref{simp},  the Cauchy-Schwarz and the Young's inequality, for any $\kappa > 1$ we have
\begin{align*}
&||\zeta_{h,t}||_a^2 + \dfrac{1}{2}\dfrac{d}{dt}||\zeta_{h}||_a^2+\dfrac{1}{2}\dfrac{d}{dt}s_h(\xi_h,\xi_h)
\nonumber \\
 =& 
a_h(({Q}_k^0\bm{\sigma - \sigma})_t, \zeta_{h,t}) +l_1(\bm{\sigma}_t, \xi_h) + l_2(\bm{v}_t, \xi_h)
 \nonumber \\
 \leq & \dfrac{1}{2}||({Q}_k^0\bm{\sigma - \sigma})_t||_a^2+\dfrac{1}{2} ||\zeta_{h,t}||_a^2 + Ch^{k+1}|\bm{\sigma}_t|_{k+1}||\bm{\varepsilon_{w,k}}(\xi_h)||_0+ Ch^{k+1}(|\bm{\sigma}_t|_{k+1}+|\bm{v}_t|_{k+2})s_h(\xi_h,\xi_h)^{1/2}
\nonumber \\
\leq & \dfrac{1}{2}||({Q}_k^0\bm{\sigma - \sigma})_t||^2 + \dfrac{1}{2}||\zeta_{h,t}||_a^2 +\dfrac{C}{2\kappa}||\bm{\varepsilon_{w,k}}(\xi_{h})||_0^2
+\dfrac{\kappa C}{2}h^{2k+2}(|\bm{\sigma}_t|_{k+1}^2+|\bm{v}_t|_{k+2}^2)+\dfrac{C}{2\kappa}s_h(\xi_h,\xi_h) 
\nonumber \\
\leq& C h^{2k+2} (|\bm{\sigma}_t|_{k+1}^2+|\bm{v}_t|_{k+2}^2) +\dfrac{1}{2}||\zeta_{h,t}||_a^2 + \dfrac{C}{2\kappa}||\bm{\varepsilon_{w,k}}(\xi_{h})||_0^2 + \dfrac{C}{2\kappa}s_h(\xi_h,\xi_h),
\end{align*}
which implies  
\begin{align} \label{335}
||\zeta_{h,t}||_a^2 +\dfrac{d}{dt}\left(||\zeta_{h}||_a^2+s_h(\xi_h,\xi_h)\right) \leq Ch^{2k+2}(|\bm{\sigma}_t|_{k+1}^2+|v_t|_{k+2}^2)+\dfrac{C}{\kappa}||\bm{\varepsilon_{w,k}}(\xi_{h})||_0^2 + \dfrac{C}{\kappa} s_h(\xi_h,\xi_h)	.
\end{align}
From \eqref{strain_estim} and  the norm equivalence \eqref{a-norm},    we have
\begin{align*}
||\bm{\varepsilon_{w,k}}(\xi_h)||_0^2   \leq&  3c^2(||\zeta_{h,t}||_0^2+||\zeta_h||_0^2+h^{2k+2}|\bm{\sigma_t}|_{k+1}^2)\\
\leq& \dfrac{3c^2}{M_0}(||\zeta_{h,t}||_a^2+||\zeta_h||_a^2)+3c^2h^{2k+2}|\bm{\sigma}_t|_{k+1}^2,
\end{align*}
which, together with \eqref{335}, yields
\begin{align*}
&||\bm{\varepsilon_{w,k}}(\xi_h)||_0^2 + \dfrac{3c^2}{M_0} \dfrac{d}{dt}\left(||\zeta_{h}||_a^2+s_h(\xi_h,\xi_h)\right) \\
\leq & \dfrac{3c^2}{M_0}\left[||\zeta_{h,t}||_a^2+\dfrac{d}{dt}(||\zeta_h||_a^2 + s_h(\xi_h,\xi_h))\right] + \dfrac{3c^2}{M_0}||\zeta_h||_a^2+3c^2h^{2k+2}|\bm{\sigma}_t|_{k+1}^2
\\
\leq & \dfrac{3c^2}{M_0}\left[ Ch^{2k+2}(|\bm{\sigma}_t|_{k+1}^2+|\bm v_t|_{k+2}^2)+\dfrac{C}{\kappa}||\bm{\varepsilon_{w,k}}(\xi_{h})||_0^2 + \dfrac{C}{\kappa} s_h(\xi_h,\xi_h)\right]+ \dfrac{3c^2}{M_0}||\zeta_h||_a^2+3c^2h^{2k+2}|\bm{\sigma}_t|_{k+1}^2.
\end{align*}
Then we get 
\begin{align} \label{336}
&\dfrac{M_0}{3c^2}(1-\dfrac{C}{\kappa})||\bm{\varepsilon_{w,k}}(\xi_h)||_0^2+  \dfrac{d}{dt}\left(||\zeta_{h}||_a^2+s_h(\xi_h,\xi_h)\right)\nonumber \\
\leq& Ch^{2k+2}(|\bm{\sigma}_t|_{k+1}^2+|\bm v_t|_{k+2}^2) + ||\zeta_h||_a^2 + \dfrac{C}{\kappa}s_h(\xi_h,\xi_h),
\end{align}
By taking a sufficiently large positive constant $\kappa$ in  this inequality and using the norm equivalence \eqref{a-norm},  from (\ref{x1}) and  (\ref{strain_estim}) it follows
\begin{align} \label{x2}
\dfrac{d}{dt}\left(||\zeta_{h}||_a^2 + s_h(\xi_h,\xi_h)\right)
 + ||\zeta_{h}||_a^2 +s_h(\xi_h, \xi_h) 
\lesssim h^{2k+2} (|\bm{\sigma}|_{k+1}^2+ |\bm{v}|_{k+2}^2+|\bm{\sigma}_t|_{k+1}^2+|\bm{v}_t|_{k+2}^2) . 
\end{align}
By the continuous Gr\"onwall's inequality (\ref{continuous_gronwall}), we can get 
\begin{align}\label{338}
&&||\zeta_{h}(t)||_a^2 + s_h(\xi_h,\xi_h) \nonumber \\  &\lesssim & ||\zeta_{h}(0)||_a^2 + s_h(\xi_h(0),\xi_h(0))+ h^{2k+2}\int_0^t(|\bm{\sigma}|_{k+1}^2  + |\bm{v}|_{k+2}^2 +|\bm{\sigma}_t|_{k+1}^2 +|\bm{v}_t|_{k+2}^2)ds.
\end{align} 
In view of \eqref{semidiscretescheme_3}, it holds $$\zeta_{h}(0)={Q}_k^0 \bm{\sigma}(0) - \bm{\sigma_{h}(0)}=0.$$ The thing left is to estimate the term $s_h(\xi_h(0),\xi_h(0))$.
To this end, we  take $w_h = \xi_h$ in   \eqref{s7} and use   Lemma \ref{approximationlemma} to get
\begin{align*}
s_h(\xi_h(0),\xi_h(0)) =& l_1(\bm{\sigma}(0), \xi_h(0)) + l_2(\bm v(0), \xi_h(0)) - b_h(\zeta_h(0),\xi_h(0))\\
 =& l_1(\bm{\sigma}(0), \xi_h(0)) + l_2(\bm v(0), \xi_h(0)) 
		\\
\lesssim& h^{k+1}|\bm{\sigma }(0)|_{k+1}\cdot ||\nabla_h \xi_h(0)||_0 + h^{k+1}|\bm{\sigma(0)}|_{k+1}||\alpha^{1/2}(Q_k^b\xi_{h0}(0)-\xi_{hb}(0))||_{\partial \mathscr{T}_h} 
\\
& +h^{k+1}|\bm{v}(0)|_{k+1}||\alpha^{1/2}(Q_k^b\xi_{h0}(0)-\xi_{hb}(0))||_{\partial \mathscr{T}_h} ,
	\end{align*}
which, together with  \eqref{s1} and  \eqref{s2}, 
leads to
\begin{align} \label{sh_estimate}
s_h(\xi_h(0),\xi_h(0)) \lesssim h^{2k+2}(|\bm{\sigma}(0)|_{k+1}^2 + |\bm v(0)|_{k+2}^2 + |\bm{\sigma}_t(0)|_{k+1}^2) .
\end{align}
Combining this estimate with \eqref{338} indicates the desired result \eqref{error_stress_est}.

Now let us prove the estimate \eqref{error_strain_est}. 
From    (\ref{s2}) and  \eqref{336} with a sufficiently large $\kappa$, we  get 
\begin{align*}
&&||\bm{\varepsilon_h}(\xi_{h0})||_0^2+\dfrac{d}{dt}(||\zeta_{h}||_a^2+s_h(\xi_h,\xi_h)) \nonumber \\
&\lesssim & ||\bm{\varepsilon_{w,k}}(\xi_h)||_0^2  +||\alpha^{1/2} (Q_k^b\xi_{h0} -\xi_{hb}) ||_{\partial\mathscr{T}_h}^2 + \dfrac{d}{dt}(||\zeta_{h}||_a^2+s_h(\xi_h,\xi_h))
\nonumber \\
&\lesssim & ||\zeta_{h}||_a^2+s_h(\xi_h,\xi_h)+ h^{2k+2}(|\bm{\sigma}_t|_{k+1}^2+|\bm{v}_t|_{k+2}^2),
\end{align*}
which, together with \eqref{error_stress_est}, yields the desired estimate for $||\bm{\varepsilon_h}(\xi_{h0})||_0^2$. This finishes the proof.
\end{proof}
\end{Lem}

Applying Lemma \ref{approximation_error},   Lemma \ref{projection_property} and the triangle inequality  gives the following  error estimate for  the semi-discrete WG scheme. 

\begin{thm}
Let  $ (\bm{\sigma}, \bm{v})\in  C^1( \bm L^2(\Omega,S)\bigcap[H^{k+1}(\Omega)]^{d\times d})\times C^1([H_0^1(\Omega)\bigcap H^{k+2}(\Omega)]^d )$ be the weak solution of system \eqref{original_formulation} and $(\bm \sigma_{h},\bm v_h)\in C^1(\Sigma_{h})\times C^1(V_h^0)$ be the solution of the WG scheme \eqref{semidiscretescheme}.
Then  
\begin{align}\label{desiredestimate2}
||\bm{\sigma - \sigma}_{h}||_0 + ||\bm{\varepsilon(v)-\varepsilon}_h (\bm{v}_{h0})||_0  
\lesssim h^{k+1}\left(\tilde M_0(\sigma, v)+\tilde M_1(\sigma, v)+\tilde M_2(\sigma, v)\right)^{1/2},  
\end{align}
where $\tilde M_0(\sigma, v),\tilde M_1(\sigma, v)$ and $\tilde M_2(\sigma, v)$ are defined in Lemma \ref{approximation_error}.
\end{thm}

\section{Fully discrete weak Galerkin method}
\subsection{Backward Euler fully discrete   scheme}
We consider a  full discretization of the quasistatic viscoelastic Maxwell model based on backward Euler scheme. 
Given a positive integer $N$, let $0=t_0<t_1<\cdots<t_N=T$ be a uniform division of time domain $[0,T]$, with $t_n=n \Delta t$ and   $\Delta t = \dfrac{T}{N}$.  For any vector or tensor-valued function $g(t)$ and any $n$, we set 
$$g^n:=g(t_n), \quad  \overline{\partial_t}g^n := \dfrac{g^n -g^{n-1}}{\Delta t}.$$

Based on the semi-discrete scheme(\ref{semidiscretescheme}),  the backward Euler
fully discrete WG scheme is given as follows: for $n= 1,2, \cdots, N$, find $(\bm{\sigma}_{h}^n, \bm{v}_h^n) = (\bm{\sigma}_{h}^n, \{\bm{v}_{h0}^n, \bm{v}_{hb}^n\}) \in \Sigma_h\times V_h^0$ such that
\begin{subequations}\label{full_discrete}
\begin{align}
a_h(\overline{\partial_t}\bm{\sigma}_{h}^n, \bm \tau_h) + a_h(\bm{\sigma}_{h}^n,\bm\tau_h) - b_h(\bm \tau_h, \bm{v}_h^n)&=0, &\forall \bm \tau_h\in \Sigma_h, \label{full_discrete1} \\
b_h(\bm{\sigma}_{h}^n, \bm w_h) + s_h(\bm{v}_h^n,\bm w_h) &= (\bm{f}^n, \bm w_{h0}), &\forall \bm w_{h}\in V_h^0, \label{full_discrete2}\\
\bm{\sigma}_{h}^0 &= {Q}_k^0\psi_0.& \label{full_discrete3}
\end{align}
\end{subequations}

\begin{thm} The fully-discrete scheme \eqref{full_discrete} has a unique solution $(\bm{\sigma}_{h}^n, \bm{v}_h^n)$   $n=1,2,\cdots,N$. 

\begin{proof}
Since this is a square system, it suffices to show the homogeneous system
\begin{align}\label{42}
\left\{\begin{array}{l}\begin{aligned}
a_h(\bm\sigma_{h}^n, \tau_h) + \Delta t a_h(\bm\sigma_{h}^n, \tau_h) - \Delta t b_h(\tau_h, \bm v_h^n) &= 0, &\forall \bm \tau_h\in \Sigma_h,
\\
b_h(\bm\sigma_{h}^n, w_h) +  s_h(\bm v_h^n, w_h) &= 0, &\forall \bm w_{h}\in V_h^0
\end{aligned}\end{array}\right.
\end{align}
 only admits  a  zero solution. In fact, taking $(\tau_h, w_h) = (\bm\sigma_{h}^n, \bm v_h^n)$ and summing up the above two equations, we  obtain
\begin{align}
(1 + \Delta t) a_h(\bm\sigma_{h}^n, \bm\sigma_{h}^n) + \Delta t s_h(\bm v_h^n, \bm v_h^n) = 0,
\end{align}
which gives $\bm\sigma_{h}^n=0 $ and $ s_h(\bm v_h^n,\bm v_h^n)=0$. Then, take $\tau_h = \bm\varepsilon_{w,k}(\bm v_h^n)$ in  the first equation of \eqref{42} leads to 
$\bm\varepsilon_{w,k}(\bm v_h^n) = 0$, which, together with $ s_h(\bm v_h^n,\bm v_h^n)=0$ and  \eqref{s2}, implies
  $ \bm v_h^n=\{\bm v_{h0}^n, \bm v_{hb}^n\}=\bm 0$. This completes the proof.
\end{proof}
\end{thm}

We have the following stability results for the fully-discrete WG 
scheme \eqref{full_discrete}. 
\begin{thm}
Assume that $\Delta t < 1$, then for any $1\leq n \leq j\leq 
N$, it holds 
\begin{align}
&\sum_{n=1}^j ||\bm\sigma_{h}^n - \bm\sigma_{h}^{n-1}||_a^2  
+   ||\bm\sigma_{h}^j||_a^2 +2\sum_{n=1}^j\Delta t ||
\bm\sigma_{h}^n||_a^2 + 2\Delta t \sum_{n=1}^j s_h(\bm 
v_h^n, \bm v_h^n)\nonumber\\
=& ||\bm \sigma_{h}^0||_a^2 + \sum_{n=1}^{j}(\bm f^n, \bm 
v_{h0}^n),\label{stable_result1} \\
&\Delta t \sum_{n=1}^j||\bm{\varepsilon_h(v_{h0}^n)}||_0^2 
+\sum_{n=1}^j s_h(\bm v_h^n - \bm v_h^{n-1}, \bm v_h^n - 
\bm v_h^{n-1})+s_h(\bm {v}_{h}^j, \bm v_h^j)\nonumber \\
 \lesssim& ||\bm{\sigma_h^{0}}||_a^2 +s_h(\bm {v}_h^{0}, 
\bm v_h^{0})+ \sum_{n=1}^j (\bm f^n, \bm v_{h0}^n)+ 
\sum_{n=1}^j (\overline{\partial _t}\bm f^n, \bm v_{h0}^n).
\label{stable_result2}
\end{align}
\begin{proof} 
Taking $(\tau_h, w_h) = (\bm\sigma_{h}^n , \bm v_h^n)$ in 
the scheme \eqref{full_discrete}, we get
\begin{align} \label{r1}
\left\{
\begin{array}{rll} a_h(\overline{\partial _t}\bm\sigma_{h}^n, 
\bm\sigma_{h}^n) + a_h(\bm\sigma_{h}^n, \bm\sigma_{h}^n)-
b_h(\bm\sigma_{h}^n,\bm v_h^n) &=& 0 ,
\\
b_h(\bm \sigma_{h}^n,\bm v_h^n) + s_h(\bm v_h^n , \bm 
v_h^n) &=& (\bm f^n, \bm v_{h0}^n).
\end{array}
\right.
\end{align}
Applying the relationship 
$2(p-q,p)=(p-q, p+q) + (p-q, p-q)$ and adding the above two 
equalities, we have
\begin{align}
\dfrac{1}{2\Delta t} ||\bm\sigma_{h}^n-\bm\sigma_{h}^{n-1}||
_a^2 + \dfrac{1}{2\Delta t}(||\bm\sigma_{h}^n||_a^2 - ||
\bm\sigma_{h}^{n-1}||_a^2 ) + ||\bm\sigma_{h}^n||_a^2 + 
s_h(\bm v_h^n, \bm v_h^n) = (\bm f^n, \bm v_{h0}^n).
\end{align}
For any $j\leq N$, summing up the above inequality with $n = 
1,2,\cdots, j$, we finally obtain the desired result 
\eqref{stable_result1}.
Applying \eqref{s2}, we get
\begin{align}\label{r2}
||\bm{\varepsilon_h (v_{h0}^n)}||_0^2 \lesssim ||
\bm{\varepsilon_{w,k}(v_{h}^n)}||_0^2 + s_h(\bm{v_h^n},
\bm{v_h^n}).
\end{align}
Using the inf-sup condition \eqref{inf-sup-dis} and the equation 
\eqref{full_discrete1}, we obtain 
\begin{align*}
||\bm{\varepsilon_{w,k}(v_h^n)}||_0 \lesssim 
&\sup_{\tau_h\in\Sigma_h}\dfrac{b_h(\tau_h,\bm{v_h^n})}{||
\tau_h||_0}=\sup_{\tau_h\in\Sigma_h}
\dfrac{ a_h(\overline{\partial _t}\bm\sigma_{h}^n, \tau_{h}) + 
a_h(\bm\sigma_{h}^n, \tau_{h})}{||\tau_h||_0} \nonumber \\
\lesssim & 
||\overline{\partial _t}\bm\sigma_{h}^n||_a + ||\bm\sigma_h^n||
_a,
\end{align*}
which, together with \eqref{r2}, yields 
\begin{align}\label{r3}
||\bm{\varepsilon_h (v_{h0}^n)}||_0^2 
\lesssim  ||\overline{\partial _t}\bm\sigma_{h}^n||_a^2 + ||
\bm\sigma_h^n||_a^2 + s_h(\bm{v_h^n},\bm{v_h^n}) .
\end{align}
In light of \eqref{r1}, we have
\begin{align*}
a_h(\overline{\partial_t}\bm\sigma_h^n,\overline{\partial_t}
\bm\sigma_h^n)+a_h(\bm\sigma_h^n,\overline{\partial_t}
\bm\sigma_h^n)-b_h(\overline{\partial_t}\bm\sigma_h^n, \bm 
v_h^n) =& 0,
\\
b_h(\overline{\partial _t}\bm\sigma_{h}^n, \bm v_h^n) + 
s_h(\overline{\partial _t}\bm v_h^n,\bm v_h^n) =& 
(\overline{\partial_t}\bm f^n, \bm v_{h0}^n).
\end{align*}
Summing up these two equalities and using the identity $2p(pq)=(p-q)^2+p^2-q^2$, we arrive at
\begin{align*}
&||\overline{\partial _t}\bm{\sigma_h^n}||_a^2 + \dfrac{1}
{2\Delta t} \left(||\bm{\sigma_h^n-\sigma_h^{n-1}}||_a^2+||
\bm{\sigma_h^n}||_a^2-||\bm{\sigma_h^{n-1}}||_a^2\right) \\
 &\qquad +\dfrac{1}{2\Delta t}\left(s_h(\bm {v}_h^n-\bm{v}
_h^{n-1},\bm {v}_h^n-\bm{v}_h^{n-1})+s_h(\bm {v}_h^n, \bm{v}
_h^n)-s_h(\bm {v}_h^{n-1}, \bm{v}_h^{n-1})
\right)=(\overline{\partial_t}\bm f^n, \bm v_{h0}^n).
\end{align*}
This identity plus \eqref{r3} implies
\begin{align*}
||\bm{\varepsilon_h}(\bm v_{h0}^n)||_0^2 
\leq& 
 C\left(
 ||\bm\sigma_h^n||_a^2+ s_h(\bm{v}_h^n,\bm{v}_h^n) -
\dfrac{1}{2\Delta t}||\bm{\sigma_h^n-\sigma_h^{n-1}}||_a^2-
\dfrac{1}{2\Delta t}(||\bm{\sigma_h^n}||_a^2-||
\bm{\sigma_h^{n-1}}||_a^2) \right. 
 \nonumber\\ 
& \left. -\dfrac{1}{2\Delta t}\left(s_h(\bm {v}_h^n - \bm {v}
_h^{n-1}, \bm v_h^n - \bm {v}_h^{n-1})+ s_h(\bm v_h^n, \bm 
v_h^n)-s_h(\bm v_h^{n-1}, \bm v_h^{n-1}) \right)\right)+
(\overline{\partial_t}\bm f^n, \bm v_{h0}^n),
\end{align*}
for $n=1,2,\cdots,j$, where $C$ is positive constant 
independent of $h, \Delta t$ and $n$. 
Thus, we have
\begin{align*}
&\Delta t\sum_{n=1}^{j} ||\bm{\varepsilon_h (v_{h0}^n)}||
_0^2+\sum_{n=1}^{j}||\bm{\sigma_h^n-\sigma_h^{n-1}}||
_a^2+||\bm{\sigma_h^j}||_a^2+\sum_{n=1}^j s_h(\bm v_h^n - 
\bm v_h^{n-1}, \bm v_h^n - \bm v_h^{n-1})+s_h(\bm {v}_h^j, 
\bm v_h^j) \\
\lesssim & \Delta t \sum_{n=1}^{j}||\bm\sigma_h^n||_a^2 +
\Delta t \sum_{n=1}^{j} s_h(\bm{v_h^n},\bm{v_h^n}) +||
\bm{\sigma_h^{0}}||_a^2 +s_h(\bm {v}_h^{0}, \bm v_h^{0}) + 
\sum_{n=1}^j (\bm f^n, \bm v_{h0}^n)+ \sum_{n=1}^j 
(\overline{\partial_t}\bm f^n, \bm v_{h0}^n).\\
\lesssim & ||\bm{\sigma_h^{0}}||_a^2 +s_h(\bm {v}_h^{0}, 
\bm v_h^{0}) + \sum_{n=1}^j (\bm f^n, \bm v_{h0}^n)+ 
\sum_{n=1}^j (\overline{\partial_t}\bm f^n, \bm v_{h0}^n),
\end{align*}
where in the second estimate we have used the stability 
result \eqref{stable_result1}. Hence, the desired result 
\eqref{stable_result2} follows.
\end{proof}
\end{thm}

\subsection{Error estimation}

By following the same line as in the proof of Lemma \ref{semi_discrete_proj_problem}, we can derive the following lemma.
\begin{Lem}\label{full_discrete_proj_lem}
Let $(\bm\sigma, \bm v)\in C^1(\bm L^2(\Omega,S)\bigcap  \bm H ({\bf div},\Omega))\times C^0([H_0^1(\Omega)]^d)$ be weak solution of system \eqref{original_formulation},  then for all $\bm{\tau_h}\in\Sigma_{h}$ and  $\bm{w_h=\{w_{h0},w_{hb}\}}\in V_h^0$, it holds
\begin{subequations} \label{full_discrete_proj_problem}
\begin{align}
&a_h(\overline{\partial_t} {Q}_k^0\bm\sigma^n ,\bm\tau_h) + a_h( {Q}_k^0\bm\sigma^n,\bm\tau_h) - b_h(\bm\tau_h,\{Q_{k+1}^0\bm v^n, Q_k^b\bm v^n\})=a_h(\overline{\partial_t} {Q}_k^0\bm\sigma^n- \bm\sigma_t^n, \bm \tau_h) \label{full_discrete_proj_problem1},
		\\
&b_h ( {Q}_k^0\bm\sigma^n,\bm w_h) + s_h(\{Q_{k+1}^0\bm v^n, Q_k^b\bm v^n\}, \bm w_h) = (\bm f^n,\bm w_h)+ l_1(\bm\sigma^n, \bm w_h) + l_2(\bm v^n,\bm w_h),\label{full_discrete_proj_problem2}
\end{align}
\end{subequations}
 for $n=1,2,\cdots,N$, where the bilinear  forms $l_1(), \ l_2()$ are defined  in Lemma \ref{semi_discrete_proj_problem}.
\end{Lem}

\begin{Lem}\label{stress_strain_est}
 Let $(\bm\sigma, \bm v)\in C^2(\bm L^2(\Omega,S)\cap [H^{k+1}(\Omega)]^{d\times d})\times C^1([H_0^1(\Omega)\cap H^{k+2}(\Omega)]^d)$ be the solution of   \eqref{original_formulation}, and let $(\bm \sigma_{h}^n, \bm v_h^n)=(\bm \sigma_{h}^n, \{\bm v_{h0}^n, \bm v_{hb}^n\})$  be the solution of   \eqref{full_discrete}  for $n=1,2,\cdots,N$. Then  
it holds
\begin{align}
||\zeta_h^n||_0^2 + 2\Delta t\sum_{j=1}^n||\zeta_h^j||_0^2 +&s_h(\xi_h^n,\xi_h^n)+2\Delta t \sum_{j=1}^n s_h(\xi_h^j,\xi_h^j)
\nonumber \\ 
\lesssim& h^{2k+2}\left(\tilde{M}_0(0)+\tilde{M}_1(t_n)+\tilde{M}_2(t_n)\right)+\Delta t^2\tilde{M}_3(t_n), \label{discrete_estimate1}
\\
\Delta t \sum_{j=1}^n||\bm\varepsilon_h(\xi_{h0}^j)||^2
\lesssim& 
h^{2k+2}\left(\tilde{M}_0(0)+\tilde{M}_1(t_n)+ \tilde{M}_2(t_n)\right)+\Delta t^2\tilde{M}_3(t_n), \label{discrete_estimate2}
\end{align}
where $\zeta_h^n : = Q_k^0\bm\sigma^n-\bm\sigma_h^n,\  \xi_h^n := \{\xi_{h0}^n,\xi_{hb}^n\}$ with $\xi_{h0}^n = Q_{k+1}^0\bm v^n-\bm v_{h0}^n$ and $\xi_{hb}^n = Q_{k}^b\bm v^n-\bm v_{hb}^n,$
\begin{align*}
\tilde{M}_0(0):=& |\bm\sigma(0)|_{k+1}^2 + |\bm v(0)|_{k+2}^2 + |\bm \sigma_t(0)|_{k+1}^2,
\\
\tilde{M}_1(t_n):=& \max_{t_j\in [0,T],1\leq j\leq n}(|\bm\sigma^j|_{k+1}^2 + |\bm v^j|_{k+2}^2),
\\
\tilde{M}_2(t_n):=& \int_0^{t_n}(|\bm\sigma_t|_{k+1}^2+|\bm v_t|_{k+2}^2)ds,
\\
\tilde{M}_3(t_n):=& \int_0^{t_n}||\bm\sigma_{tt}||_{0}^2ds.
\end{align*}
\begin{proof} The proof is similar to that of Lemma \ref{approximation_error} for   the semi-discrete scheme.  For completeness, we show it as following.
We mention that the notation $C_i$  in this proof for any $i$ denotes a generic  positive constant independent of $h$ and $ \Delta t$. 

Our proof mainly divides into 4 steps.

{$\bm Step $} 1: 
From (\ref{full_discrete}) and (\ref{full_discrete_proj_problem}) it follows, for any  $1\leq j\leq n$,  
\begin{subequations}
\begin{align}
a_h(\overline{\partial _t}\zeta_{h}^j, \tau_h) + a_h(\zeta_{h}^j,\tau_h) - b_h(\tau_h, \xi_h ^j) &= a_h(\overline{\partial_t} {Q}_{k}^0\bm\sigma^j- \bm\sigma_t^j, \tau_h), \label{a1} \\
b_h(\zeta_{h}^j,w_h) + s_h(\xi_h^j, w_h) &= l_1(\bm\sigma^j,w_h) + l_2(\bm v^j, w_h). \label{a2}
\end{align}
\end{subequations}
Taking $\tau_h=\zeta_h^j$ and $w_h=  \xi_h^j$, and summing up the above two equations, we obtain
\begin{align}\label{discrete_error}
a_h(\overline{\partial _t}\zeta_{h}^j, \zeta_{h}^j) + a_h(\zeta_{h}^j,\zeta_{h}^j) + s_h(\xi_h^j, \xi_h ^j) &= 
 l_1(\bm\sigma^j,\xi_h ^j) 
+ l_2(\bm v^j, \xi_h^j) + a_h(\overline{\partial_t} {Q}_k^0\bm{\sigma}^j- \bm\sigma_t^j, \zeta_{h}^j)\nonumber\\
&=:E_1^j + E_2^j + E_3^j.
\end{align}

Taking $\tau_h = \overline{\partial _t} \zeta_h^j$ in equality \eqref{a1} to get
\begin{align*}
	a_h(\overline{\partial_t}\zeta_h^j,\overline{\partial_t}\zeta_h^j)+a_h(\zeta_h^j,\overline{\partial_t}\zeta_h^j)-b_h(\overline{\partial_t}\zeta_h^j, \xi_h^j) = a_h(\overline{\partial_t}(Q_k^0\bm\sigma^j-\bm\sigma^j), \overline{\partial_t}\zeta_h^j).
\end{align*}
In light of   \eqref{full_discrete_proj_problem2} and the fact that $-\nabla\cdot\overline{\partial _t}\bm\sigma^n = \overline{\partial_t}f^n$, we  have
\begin{align*}
b_h(\overline{\partial _t}\zeta_h^j,\xi_h^j) + s_h(\overline{\partial _t}\xi_h^j,\xi_h^j) = l_1(\overline{\partial _t}\bm\sigma^j,\xi_h^j) + l_2(\overline{\partial _t}\bm v^j, \xi_h^j).
\end{align*}
Summing up the above two equalities, we  obtain
\begin{align}\label{discrete_error_t}
	a_h(\overline{\partial_t}\zeta_h^j,\overline{\partial_t}\zeta_h^j)+a_h(\zeta_h^j,\overline{\partial_t}\zeta_h^j)  + s_h(\overline{\partial _t}\xi_h^j,\xi_h^j) 
	=
	a_h(\overline{\partial_t}({Q}_k^0\bm\sigma^j-\bm\sigma^j), \overline{\partial_t}\zeta_h^j)+ l_1(\overline{\partial _t}\bm\sigma^j,\xi_h^j) + l_2(\overline{\partial _t}\bm v^j, \xi_h^j), 
\end{align}
which shows 
 \begin{align*}
& ||\overline{\partial_t}\zeta_h^j||_a^2+a_h(\zeta_h^j,\overline{\partial_t}\zeta_h^j)+s_h(\overline{\partial_t}\xi_h^j,\xi_h^j)\nonumber \\
 \leq &
       \dfrac{1}{2}||\overline{\partial_t}(Q_k^0\bm\sigma^j-\bm\sigma^j)||_a^2 + \dfrac{1}{2}||\overline{\partial_t}\zeta_h^j||_a^2  +||\alpha^{1/2}(Q_k^b\xi_{h0}^j - \xi_{h0}^j)||_{\partial \mathscr{T}_h}||\alpha^{-1/2}\overline{\partial_t}( Q_k^0\bm\sigma^j - \bm\sigma^j)  ||_{\partial \mathscr{T}_h}
       \nonumber \\
& 
      +||\alpha^{1/2}(Q_k^b\xi_{h0}^j - \xi_{hb}^j)||_{\partial \mathscr{T}_h}\left(||\alpha^{-1/2}\overline{\partial_t}(Q_k^0\bm\sigma^j -  \bm\sigma^j)  ||_{\partial \mathscr{T}_h}
      + ||\alpha^{1/2}Q_k^b\overline{\partial _t}(Q_{k+1}^0\bm v^j-\bm v^j)||_{\partial\mathscr{T}_h}\right).
 \end{align*}
Thus, 
 \begin{align}\label{417}
 & \frac12||\overline{\partial_t}\zeta_h^j||_a^2+a_h(\zeta_h^j,\overline{\partial_t}\zeta_h^j)+s_h(\overline{\partial_t}\xi_h^j,\xi_h^j)\nonumber \\
\leq &
       \dfrac{1}{2}||\overline{\partial_t}(Q_k^0\bm\sigma^j-\bm\sigma^j)||_a^2 + ||\alpha^{1/2}(Q_k^b\xi_{h0}^j - \xi_{h0}^j)||_{\partial \mathscr{T}_h}||\alpha^{-1/2}\overline{\partial_t}( Q_k^0\bm\sigma^j - \bm\sigma^j)  ||_{\partial \mathscr{T}_h}
       \nonumber \\
& 
      +||\alpha^{1/2}(Q_k^b\xi_{h0}^j - \xi_{hb}^j)||_{\partial \mathscr{T}_h}\left(||\alpha^{-1/2}\overline{\partial_t}( Q_k^0\bm\sigma^j - \bm\sigma^j)  ||_{\partial \mathscr{T}_h}
      + ||\alpha^{1/2}Q_k^b\overline{\partial _t}(Q_{k+1}^0\bm v^j-\bm v^j)||_{\partial\mathscr{T}_h}\right).
 \end{align}
By the projection properties of ${Q}_k^0$, we obtain
\begin{align}\label{partial_sigma_error}
||\overline{\partial_t}{Q}_k^0\bm\sigma^j - \overline{\partial_t} \bm\sigma^j||_{0,\mathscr{T}_h} =& ||{Q}_k^0\overline{\partial_t}\bm\sigma^j - \overline{\partial_t} \bm\sigma^j||_{0,\mathscr{T}_h} = \dfrac{1}{\Delta t} {\int}_{t_{j-1}}^{t_j} |{Q}_k^0\bm\sigma_t(s)-\bm\sigma_t(s)|ds \nonumber \\
\lesssim & \dfrac{h^{k+1}}{\Delta t} {\int}_{t_{j-1}}^{t_j}|\bm\sigma_t(s)|_{k+1}ds
\lesssim  \dfrac{h^{k+1}}{\sqrt{\Delta t}}\left(\int_{t_{j-1}}^{t_j}|\bm\sigma_t(s)|_{k+1}^2ds\right)^{\frac{1}{2}},\\
||\overline{\partial_t}{Q}_k^0\bm\sigma^j - \overline{\partial_t} \bm\sigma^j||_{\partial\mathscr{T}_h} 
\lesssim&  \dfrac{h^{k+\frac12}}{\sqrt{\Delta t}} \left(\int_{t_{j-1}}^{t_j}|\bm\sigma_t(s)|_{k+1}^2ds \right)^{\frac{1}{2}}.
\end{align}
Similarly, we have 
\begin{align}\label{partial_v_error-bd}
||\overline{\partial_t}{Q}_k^0\bm v^j - \overline{\partial_t} \bm v^j||_{\partial\mathscr{T}_h} 
\lesssim  \dfrac{h^{k+\frac32}}{\sqrt{\Delta t}}\left(\int_{t_{j-1}}^{t_j}|\bm v_t(s)|_{k+2}^2ds\right)^{\frac{1}{2}}.
\end{align}
In light of  \eqref{partial_sigma_error}-\eqref{partial_v_error-bd}, the the projection properties of $Q_k^b$, the Young's inequality,    
  the norm equivalence \eqref{a-norm} and  Lemma \ref{simp}, we further apply \eqref{417} to get
\begin{align}
&\frac12||\overline{\partial _t}\zeta_h^j||_a^2 + a_h(\zeta_h^j,\overline{\partial _t}\zeta_h^j) + s_h(\overline{\partial _t}\xi_h^j,\xi_h^j) 
\nonumber\\
\leq &  \dfrac{1}{p} ||\bm{\varepsilon_{w,k}}(\xi_h^j)||_0^2 + \dfrac{1}{p}s_h(\xi_h^j, \xi_h^j) + C_1\dfrac{h^{2k+2}}{\Delta t}\int_{t_{j-1}}^{t_j}\left( |\bm\sigma_t(s)|_{k+1}^2+|\bm v_t(s)|_{k+2}^2\right)ds,
\end{align}
where $p>0$ is an arbitrary positive constant. 

A similar proof of  $\bm{\varepsilon_{w,k}}(\xi_h^j)$ as that of \eqref{strain_estim} implies
\begin{align} \label{discrete_strain1}
||\bm{\varepsilon_{w,k}}(\xi_{h}^j)||_0^2 
\leq C_2\left(||\overline{\partial _t}\zeta_h^j||_a^2+||\zeta_h^j||_a^2  + \dfrac{h^{2k+2}}{\sqrt{\Delta t}}\int_{t_{j-1}}^{t_j}|\bm\sigma_t|_{k+1}^2\right) . 
\end{align}
Hence, if we choose $p$ sufficiently large such that $p>4C_2$, then the above two inequalities give
\begin{align} \label{zeta_hj_estimate}
&\frac14||\overline{\partial _t}\zeta_h^j||_a^2 + a_h(\zeta_h^j,\overline{\partial _t}\zeta_h^j) + s_h(\overline{\partial _t}\xi_h^j,\xi_h^j) 
\nonumber\\
\leq &   \dfrac{1}{p}s_h(\xi_h^j, \xi_h^j) + \dfrac{C_2}{p}||\zeta_h^j||_a^2 +(C_1+ \dfrac{C_2}{p})\dfrac{h^{2k+2}}{\Delta t}\int_{t_{j-1}}^{t_j}\left( |\bm\sigma_t(s)|_{k+1}^2+|\bm v_t(s)|_{k+2}^2\right),
\end{align}
and from \eqref{discrete_strain1} we get
\begin{align}\label{discrete_strain2}
||\bm{\varepsilon_{w,k}}(\xi_{h}^j)||_0^2  \leq& 
C_3\left(s_h(\xi_h^j,\xi_h^j) +  ||\zeta_h^j||_a^2 +  \dfrac{h^{2k+2}}{\Delta t}\int_{t_{j-1}}^{t_j}(|\bm\sigma_t(s)|_{k+1}^2+|\bm v_t(s)|_{k+2}^2)ds \right.\nonumber \\
&\qquad \left. \quad -a_h(\zeta_h^j, \overline{\partial _t}\zeta_h^j)-s_h(\overline{\partial _t}\xi_h^j,\xi_h^j)\right).
\end{align}

{$\bm Step$ 2}: The next thing is to estimate the terms $E_1^j, E_2^j$ and $E_3^j$ in \eqref{discrete_error}, respectively.  From Lemma \ref{approximationlemma} and Lemma \ref{simp},  it follows
\begin{align*} 
E_1^j  =l_1(\bm\sigma^j,\xi_h^j)     \lesssim& h^{k+1}|\bm\sigma^j|_{k+1}\left(||\bm{\varepsilon_{w,k}}(\xi_{h}^j)||_0 +s_h(\xi_h^j,\xi_h^j)^{1/2} \right),
\end{align*}
which, together with the 
 Cauchy inequality, the  Young's inequality and \eqref{discrete_strain2}, indicates
\begin{align} \label{E1}
E_1^j     
\leq  &    C_4 \left(h^{2k+2}|\bm\sigma^j|_{k+1}^2 +  \dfrac{h^{2k+2}}{\Delta t} \int_{t_{j-1}}^{t_j}(|\bm\sigma_t|_{k+1}^2  + |\bm v_t|_{k+2}^2)ds -a_h(\overline{\partial_t}\zeta_{h}^j,\zeta_{h}^j) - s_h(\overline{\partial _t}\xi_h^j,\xi_h^j)\right)
\nonumber \\
     & \quad + \dfrac{1}{p'}||\zeta_h^j||_a^2  + \dfrac{1}{p'}s_h(\xi_h^j,\xi_h^j) .
\end{align}
where $p'$ is an arbitrary positive constant.
For the term $E_2^j$,  by Lemma \ref{approximationlemma}  and  the  Young's inequality we also have
\begin{align}\label{E2}
E_2^j = l_2(\bm v^j, \xi_h^j)  \leq  C_5 h^{2k+2} |\bm v^j|_{k+2}^2 + \dfrac{1}{p'}s_h(\xi_h^j,\xi_h^j).
\end{align}
Applying the Taylor formula and the Cauchy inequality  gives
\begin{align*}
||\overline{\partial_t} \sigma^j-\sigma_t^j||_{0,\mathscr{T}_h} = \dfrac{1}{\Delta t} {\int}_{t_{j-1}}^{t_j} (s-t_{j-1})||\sigma_{tt}(s)||_0ds  \leq  \sqrt{\Delta t} \left({\int}_{t_{j-1}}^{t_j} ||\bm\sigma_{tt}(s)||_0^2 ds \right )^{\frac{1}{2}},
\end{align*}
which, together with \eqref{partial_sigma_error}, shows
\begin{align}\label{E3}
E_3^j =& a_h(\overline{\partial_t} {Q}_k^0\bm\sigma^j- \bm\sigma_t^j, \zeta_{h}^j) = a_h(\overline{\partial_t} {Q}_k^0\bm\sigma^j - \overline{\partial_t} \bm\sigma^j, \zeta_{h}^j) + a_h(\overline{\partial_t} \bm\sigma^j-\bm\sigma_t^j, \zeta_{h}^j) \nonumber \\
 \leq & M_1\left(||\overline{\partial_t}{Q}_k^0\bm\sigma^j - \overline{\partial_t} \bm\sigma^j||_{0,\mathscr{T}_h}+||\overline{\partial_t} \bm\sigma^j-\bm\sigma_t^j||_{0,\mathscr{T}_h}\right)||\zeta_{h}^j||_{0,\mathscr{T}_h} \nonumber\\
\leq & C_6\left( \dfrac{h^{2k+2}}{ \Delta t} \int_{t_{j-1}}^{t_j}|\bm\sigma_t|_{k+1}^2ds  +  {\Delta t}  {\int}_{t_{j-1}}^{t_j}||\bm\sigma_{tt}||_0^2 ds  \right)+\frac{1}{p'}||\zeta_{h}^j||_{a}^2.
\end{align}

{$\bm Step $ 3:} The equality \eqref{discrete_error} plus the estimates   \eqref{E1}-\eqref{E3}   implies  
\begin{align*}
&a_h(\overline{\partial _t}\zeta_h^j,\zeta_h^j) +a_h(\zeta_h^j,\zeta_h^j) +s_h(\xi_h^j,\xi_h^j)
\nonumber \\
\leq & 
C_7\left( h^{2k+2}(|\bm\sigma^j|_{k+1}^2+|\bm v^j|_{k+2}^2) + \dfrac{h^{2k+2}}{\Delta t} \int_{t_{j-1}}^{t_j}(|\bm\sigma_t|_{k+1}^2  + |\bm v_t|_{k+2}^2)ds+   {\Delta t}  {\int}_{t_{j-1}}^{t_j}||\bm\sigma_{tt}||_0^2 ds \right)\nonumber \\
& \quad    + \dfrac{2}{p'}||\zeta_h^j||_a^2  + \dfrac{2}{p'}s_h(\xi_h^j,\xi_h^j)  - C_4\left(a_h(\overline{\partial_t}\zeta_{h}^j,\zeta_{h}^j) +s_h(\overline{\partial _t}\xi_h^j,\xi_h^j)\right).
\end{align*}
Taking $p'=4$ in this inequality, we further obtain 
\begin{align*}
&(1+C_4)a_h(\overline{\partial _t}\zeta_h^j,\zeta_h^j) +C_4s_h(\overline{\partial _t}\xi_h^j,\xi_h^j)+\frac12||\zeta_h^j||_a^2 +\frac12s_h(\xi_h^j,\xi_h^j)
\nonumber \\
\leq & 
C_7\left( h^{2k+2}(|\bm\sigma^j|_{k+1}^2+|\bm v^j|_{k+2}^2) + \dfrac{h^{2k+2}}{\Delta t} \int_{t_{j-1}}^{t_j}(|\bm\sigma_t|_{k+1}^2  + |\bm v_t|_{k+2}^2)ds +  {\Delta t}  {\int}_{t_{j-1}}^{t_j}||\bm\sigma_{tt}||_0^2 ds\right) ,
\end{align*}
which means
\begin{align*}
& a_h(\overline{\partial _t}\zeta_h^j,\zeta_h^j) + s_h(\overline{\partial _t}\xi_h^j,\xi_h^j)+ ||\zeta_h^j||_a^2 + s_h(\xi_h^j,\xi_h^j)
\nonumber \\
\leq & 
C_8\left( h^{2k+2}(|\bm\sigma^j|_{k+1}^2+|\bm v^j|_{k+2}^2) + \dfrac{h^{2k+2}}{\Delta t} \int_{t_{j-1}}^{t_j}(|\bm\sigma_t|_{k+1}^2  + |\bm v_t|_{k+2}^2)ds +  {\Delta t}  {\int}_{t_{j-1}}^{t_j}||\bm\sigma_{tt}||_0^2 ds  \right).
\end{align*}
This estimate, together with the identities
\begin{align*}
	&a_h(\overline{\partial_t} \zeta_{h}^j, \zeta_{h}^j) =\dfrac{1}{2\Delta t} \left(||\zeta_{h}^j-\zeta_{h}^{j-1}||_a^2 + ||\zeta_{h}^j||_a^2-||\zeta_{h}^{j-1}||_a^2\right), \\
	&s_h(\overline{\partial _t}\xi_h^j,\xi_h^j)=\dfrac{1}{2\Delta t}(s_h(\xi_h^j-\xi_h^{j-1}, \xi_h^j-\xi_h^{j-1}) + s_h(\xi_h^j, \xi_h^j)-s_h(\xi_h^{j-1},\xi_h^{j-1})),
\end{align*}
yields
\begin{align*}
&||\zeta_h^j||_a^2-||\zeta_{h}^{j-1}||_a^2 +  s_h(\xi_h^j,\xi_h^j) - s_h(\xi_h^{j-1}, \xi_h^{j-1})  + 2\Delta t  ||\zeta_{h}^j||_a^2 + 2\Delta t  s_h(\xi_h^j,\xi_h^j)
\nonumber \\
\leq &
  2 C_8\left( \Delta t h^{2k+2}(|\bm\sigma^j|_{k+1}^2+|\bm v^j|_{k+2}^2) + h^{2k+2}  \int_{t_{j-1}}^{t_j}(|\bm\sigma_t|_{k+1}^2  + |\bm v_t|_{k+2}^2)ds  +  (\Delta t) ^2 {\int}_{t_{j-1}}^{t_j}||\bm\sigma_{tt}||_0^2 ds  \right).
\end{align*}
Summing up the above inequality for $j=1,\cdots,n$, we arrive at
\begin{align}
& ||\zeta_h^n||_a^2 +  s_h(\xi_h^n,\xi_h^n)+{2\Delta t}\sum_{j=1}^n\left( ||\zeta_{h}^j||_a^2+  s_h(\xi_h^j,\xi_h^j) \right)
\nonumber \\ 
\leq&
      ||\zeta_{h}^0||_a^2 +  s_h(\xi_h^0,\xi_h^0) + 2C_8\left(h^{2k+2} \max_{t_j\in[0,T]}\left( |\bm\sigma^j|_{k+1}^2+|\bm v^j|_{k+2}^2  \right)\right.\nonumber\\
      &\quad \left.+h^{2k+2}\int_0^{t_n}(|\bm\sigma_t|_{k+1}^2+|\bm v_t|_{k+2}^2)ds +(\Delta t)^2 \int_0^{t_n}||\bm\sigma_{tt}||_0^2ds\right),
\end{align}
which, together with \eqref{full_discrete3} and \eqref{sh_estimate}, leads to desired estimate \eqref{discrete_estimate1}. 

{$\bm Step $ 4: }
Finally, let us prove  \eqref{discrete_estimate2}.
From inequality (\ref{s2}) and \eqref{discrete_strain2},  we get
\begin{align*}
||\bm{\varepsilon_h}(\xi_{h0}^j)||_0^2 \lesssim& ||\bm{\varepsilon_{w,k}}(\xi_{h}^j)||_0^2 + s_h(\xi_h^j,\xi_h^j)
\nonumber \\
\lesssim & 
C_3 \left( s_h(\xi_h^j,\xi_h^j) + ||\zeta_h^j||_a^2+\dfrac{h^{2k+2}}{\Delta t} \int_{t_{j-1}}^{t_j} (|\bm\sigma_t|_{k+1}^2+|\bm v_t|_{k+2}^2)ds
\right. \nonumber \\
&\left.-a_h(\overline{\partial_t}\zeta_{h}^j, \zeta_{h}^j)  - s_h(\overline{\partial_t}\xi_h^j,\xi_h^j)\right)+ s_h(\xi_h^j,\xi_h^j),
\end{align*}
which implies 
\begin{align*}
&||\varepsilon_h(\xi_{h0}^j)||_0^2 + \dfrac{1}{2\Delta t}\left (||\zeta_h^j||_a^2-||\zeta_h^{j-1}||_a^2+s_h(\xi_h^j,\xi_h^j)-s_h(\xi_h^{j-1},\xi_h^{j-1})\right) \nonumber \\
\leq&  C_9 \left(\dfrac{h^{2k+2}}{\Delta t} \int_{t_{j-1}}^{t_j} (|\bm\sigma_t|_{k+1}^2+|\bm v_t|_{k+2}^2)ds
+ s_h(\xi_h^j,\xi_h^j) + ||\zeta_h^j||_a^2\right) .
\end{align*}
Summing up the above inequality for $j=1,2,\cdots,n$, we have
\begin{align}
&2\Delta t \sum_{j=1}^n||\bm{\varepsilon_h}(\xi_{h0}^j)||_0^2
 + ||\zeta_{h}^n||_a^2 
 + s_h(\xi_h^n,\xi_h^n)
\nonumber \\
\lesssim & ||\zeta_h^0||_a^2+ s_h(\xi_h^0,\xi_h^0)+ h^{2k+2}\int_0^{t_n}(|\bm\sigma_t|_{k+1}^2+|\bm v_t|_{k+2}^2)ds
 +  2\Delta t\sum_{j=1}^n (||\zeta_{h}^j||_a^2 + s_h(\xi_h^j,\xi_h^j)),
\end{align}
which,  together  with \eqref{full_discrete3}, \eqref{sh_estimate} and  \eqref{discrete_estimate1}, yields
the desired result (\ref{discrete_estimate2}).
\end{proof}
\end{Lem}

Applying Lemmas \ref{stress_strain_est}, Lemma \ref{projection_property} and the triangle inequality leads to  the following error estimate for the  fully discrete scheme.
\begin{thm}\label{fully-convergence}
Let $(\bm\sigma, \bm v)\in C^2(\bm L^2(\Omega,S)\cap [H^{k+1}(\Omega)]^{d\times d})\times C^1([H_0^1(\Omega)\cap H^{k+2}(\Omega)]^d)$ be the solution of   \eqref{original_formulation}, and let $(\bm \sigma_{h}^n, \bm v_h^n)=(\bm \sigma_{h}^n, \{\bm v_{h0}^n, \bm v_{hb}^n\})$  be the solution of   \eqref{full_discrete}  for $n=1,2,\cdots,N$. Then 
it holds
\begin{align}
||\bm\sigma(t_n)-\bm\sigma_{h}^n||_{0}^2 + \Delta t ||\bm{\varepsilon (v(t_n) ) }- \bm{\varepsilon_h( v_{h0}^n)}||^2 \lesssim h^{2k+2}\left(\tilde{M}_0(0)+\tilde{M}_1(t_n)+\tilde{M}_2(t_n)\right)+\Delta t^2\tilde{M}_3(t_n),
\end{align}
where $\tilde{M}_0(0),\tilde{M}_1(t_n),\tilde{M}_2(t_n)$ and $\tilde{M}_3(t_n)$ are defined in Lemma \ref{stress_strain_est}.
\end{thm}

\section{Numerical examples}

In this section, we provide two  2-dimensional numerical examples  to verify the performance of the proposed fully discrete weak Galerkin method \eqref{full_discrete} with $k=1,2$. 

In the numerical examples of  the  model problem \eqref{modelproblem},  we  take   $\Omega = [0,1]\times[0,1]$ and $T =1 $, and assume  the elastic medium to be isotropic with   $\mu=1 $ and $\lambda = 1$. We use $M\times M$ uniform triangular meshes (c.f. Figure 1) for the spatial discretization.

\centerline{\includegraphics[height=6cm,width=12cm]{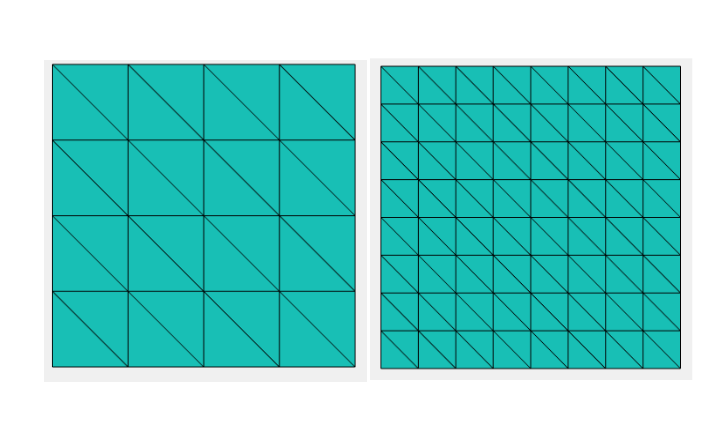}}
\centerline{Figiure 1: The domain:  $4\times4$ (left)  and $8\times 8$ (right) meshes}

\begin{Exa}\label{Exa5.1} 
The exact displacement field $\bm{u}(x,t)$ and symmetric stress tensor $\bm{\sigma}(x,t) = (\sigma_{ij})_{2\times2}$ are respectively given by
$$\bm{u}=\left(\begin{array}{cc}-e^{-t}(x_1^4-2x_1^3+x_1^2)(4x_2^3-6x_2^2+2x_2) \\-e^{-t}(x_2^4-2x_2^3+x_2^2)(4x_1^3-6x_1^2+2x_1)\end{array}\right),$$
	$$\left(\begin{array}{cc} \sigma_{11} \\ \sigma_{12} \\ \sigma_{22} \end{array} \right) = \left(\begin{array}{cc}16te^{-t}(2x_1^3-3x_1^2+x_1)(2x_2^3-3x_2^2+x_2)\\ 2te^{-t}[(x_1^4-2x_1^3+x_1^2)(6x_2^2-6x_2+1)+(x_2^4-2x_2^3+x_2^2)(6x_1^2-6x_1+1)]\\ 16te^{-t}(2x_1^3-3x_1^2+x_1)(2x_2^3-3x_2^2+x_2)\end{array}\right). $$
Notice that the velocity field $\bm v=\bm u_t$.  	

To verify the spatial accuracy, we take the time step $\Delta t = 0.0005$ for $k=1$ and $\Delta t= 0.00005$ for $k=2$, respectively. Numerical results of relative errors for the discrete stress $\bm \sigma_h$ and discrete strain $\bm\varepsilon_h(\bm v_h)$ at the final time $t=T$ are presented in the Tables \ref{table1} and \ref{table2}. We can see that spatial convergence orders of the stress and strain are 
  $(k+1)-th$, which is conformable to the theoretical prediction in Theorem \ref{fully-convergence}. 

To test the temporal accuracy, we  use a very fine spatial mesh with $M=64$.  Numerical results of the errors at the final time $T=1$ are listed in Table \ref{table3}. 
We can observe  the  first order temporal  convergence rate for the stress approximation, as is consistent with Theorem    \ref{fully-convergence}, and a better rate than first order for the strain approximation.
\end{Exa}

\begin{table}[htbp]
	\centering
	\caption{History of convergence for Example \ref{Exa5.1} with   $\Delta t = 0.0005$: spatial accuracy }
	\label{table1}
\begin{tabular}{|c|c|c|c|c|c|}
	\hline
	&mesh    & $\frac{\lVert \sigma(T)- \sigma_{h}(T)\rVert_0}{\lVert \sigma(T)\rVert_0}$ & order&$\frac{\sqrt{\Delta t}||\varepsilon(v(T))-\varepsilon_h(v_h(T))||_0}{||\varepsilon(v(T))||_0}$&order
	\\ 	\hline 
	\multirow{5}{*}{$k=1$} & $ 2\times 2 $	 &1.2181e-01	    &--  &1.3172e-02&--  
	\\  \cline{2-6}
	
	& $ 4\times 4 $	 &3.3882e-02	    &1.85  &3.4212e-03& 1.94 
	\\  \cline{2-6}
	&$ 8\times 8 $	        &8.7967e-03       &1.95       & 8.7817e-04    &   1.96
	\\ \cline{2-6}
	&$16\times 16$	    &2.2206e-03      &1.99	  &2.2718e-04   &1.95
	\\ \cline{2-6}
	&$32\times 32$	    &5.5614e-04	     &2.00	&6.0670e-05 &1.90 \\
		\hline
	\end{tabular}
\end{table}

\begin{table}[htbp]
	\centering
	\caption{History of convergence for Example \ref{Exa5.1} with   $\Delta t = 0.00005$: spatial accuracy}
	\label{table2}
\begin{tabular}{|c|c|c|c|c|c|}
	\hline
	& mesh     & $\frac{\lVert \sigma(T)- \sigma_{h}(T)\rVert_0}{\lVert \sigma(T)\rVert_0}$ & order&$\frac{\sqrt{\Delta t}||\varepsilon(v(T))-\varepsilon_h(v_h(T))||_0}{||\varepsilon(v(T))||_0}$&order
	\\ 	\hline 
	\multirow{5}{*}{$k=2$} 	& $ 2\times 2 $	 &  1.3765e-01	 &--  &7.4698e-03&--  
	\\  \cline{2-6}
	
	& $ 4\times 4 $	 &  3.0684e-02	    & 2.17 &  1.3196e-03  &   2.50  
	\\  \cline{2-6}
	&$ 8\times 8 $	        &4.3824e-03       &2.81       & 1.7993e-04    &   2.87
	\\ \cline{2-6}
	&$16\times 16$	    &5.6970e-04      &2.94	  &2.3084e-05   &2.96
	\\ \cline{2-6}
	&$32\times 32$	    &7,2100e-05	     &2.98	&2.9258e-06 &2.98 \\
	\hline
\end{tabular}
\end{table}

\begin{table}[htbp] 
	\centering
	\caption{History of convergence for Example \ref{Exa5.1} with $M=64$: temporal accuracy}
	\label{table3}
		\begin{tabular}{|c|c|c|c|c|c|c|c|}
			\hline
			\multirow{8}{*}{$k=1$} 	& $\Delta t$    & $\frac{\lVert v(T)-v_{h0}(T) \rVert_0}{\lVert v(T)\rVert_0}$ & order &$\frac{\sqrt{\Delta t}||\varepsilon(v(T))-\varepsilon_h(v_h(T))||_0}{||\varepsilon(v(T))||_0}$&order 
			\\ 	\hline
			&$ 0.5 $   &  3.5128e-01   &--   &2.4956e-01 & --          
			\\
			 \cline{2-6}
			
			&$ 0.25 $   & 1.4792e-01    &  1.25 &7.4798e-02 &1.74          
			\\
			 \cline{2-6}
			&$ 0.125 $	&   6.7960e-02   &1.12     &2.4631e-02& 1.60  
			\\
			 \cline{2-6}
			&$0.0625$	&  3.2583e-02    &1.06      &8.5896e-03 &  1.52
			\\
			 \cline{2-6}
			&$0.03125$	&  1.5954e-02   &1.03    &3.1565e-03 &    1.44
			\\
			 \cline{2-6}
			&$0.015625$	&  7.8951e-03  &1.01     &1.2511e-03 &   1.34
			\\
			\hline
	\multirow{6}{*}{$k=2$}	&$ 0.5 $   &  3.5128e-01    &--   &2.4839e-01 & --          
\\
 \cline{2-6}

&$ 0.25 $   &  1.4792e-01    &1.25   &7.3962e-02 &1.75          
	\\
	 \cline{2-6}
&$ 0.125 $	&   6.7961e-02   &1.12     &2.4028e-02& 1.62  
			\\
			 \cline{2-6}
&$0.0625$	&  3.2583e-02    &1.06      &8.1460e-03 &  1.56
			\\
			 \cline{2-6}
&$0.03125$	&  1.5954e-02   &1.03    &2.8205e-03 &    1.53
			\\
			 \cline{2-6}
&$0.015625$	&  7.8944e-03  &1.02     &9.8692e-04 &   1.52
			\\
			\hline
	
	\end{tabular}
\end{table}

\begin{Exa}\label{Exa5.2} The exact displacement field $u$ and symmetric stress tensor $\sigma$ are of the following forms:
	
$$\bm{u}=\left(\begin{array}{cc}
-e^{-t}\sin (\pi x_1)\sin(\pi x_2) \\-e^{-t}\sin(\pi x_1)\sin(\pi x_2)\end{array}\right)$$
	
$$\left(\begin{array}{cc} \sigma_{11} \\ \sigma_{12} \\ \sigma_{22} \end{array} \right) = \left(\begin{array}{cc}\pi te^{-t}(3\cos(\pi x_1)\sin(\pi x_2) + \sin(\pi x_1)\cos(\pi x_2))\\ \pi te^{-t}(\sin(\pi x_1)\cos(\pi x_2) + \cos(\pi x_1)\sin(\pi x_2))\\ \pi te^{-t}(3\sin(\pi x_1)\cos(\pi x_2) + \cos(\pi x_1)\sin(\pi x_2))\end{array}\right). $$
 Tables \ref{table4}-\ref{table5} show that the scheme (\ref{full_discrete}) yields  the $(k+1)-th$ spatial convergence orders for the  the stress and strain approximations, and  Table \ref{table6} shows the  first order temporal  convergence rate for the stress approximation. These are conformable to Theorem \ref{fully-convergence}. In particular,  Table \ref{table6} also shows  a better convergence rate than first order for the strain approximation.

\end{Exa}
\begin{table}[htbp]
	\centering
	\caption{History of convergence for Example \ref{Exa5.2} with   $\Delta t = 0.0005$: spatial accuracy }
	\label{table4}
	\begin{tabular}{|c|c|c|c|c|c|}
		\hline
		& mesh     & $\frac{\lVert \sigma(T)- \sigma_{h}(T)\rVert_0}{\lVert \sigma(T)\rVert_0}$ & order&$\frac{\sqrt{\Delta t}||\varepsilon(v(T))-\varepsilon_h(v_h(T))||_0}{||\varepsilon(v(T))||_0}$&order
		\\ 	\hline 
	\multirow{4}{*}{$k=1$}  
&$ 2\times 2 $    &1.2181e-01    &--  & 1.3172e-02 &--  \\
\cline{2-6}
&$ 4\times 4 $    &   3.3882e-02 &1.85  &3.4212e-03&1.95  \\
\cline{2-6}
&$ 8\times 8 $	      &8.7967e-03    &1.95  & 8.7817e-04& 1.96\\
\cline{2-6}
&$ 16\times 16 $	    &2.2206e-03    &1.99   &2.2718e-04 & 1.96\\
\cline{2-6}
&$32\times 32$	 & 5.5614e-04   &2.00 &6.0607e-05&1.91
\\   \hline
	\end{tabular}
\end{table}

\begin{table}[htbp]
	\centering
	\caption{History of convergence for Example \ref{Exa5.2} with  $\Delta t = 0.00005$: spatial accuracy}
	\label{table5}
\begin{tabular}{|c|c|c|c|c|c|}
\hline
&mesh     & $\frac{\lVert \sigma(T)- \sigma_{h}(T)\rVert_0}{\lVert \sigma(T)\rVert_0}$ & order&$\frac{\sqrt{\Delta t}||\varepsilon(v(T))-\varepsilon_h(v_h(T))||_0}{||\varepsilon(v(T))||_0}$&order
\\ 	\hline 
\multirow{5}{*}{$k=2$} &$ 2\times 2 $    & 24575e-02   &--  & 1.1239e-03  &--  \\
\cline{2-6}
& $ 4\times 4 $	       &	3.3115e-03    & 2.89  &1.4900e-04 &2.91 
\\  \cline{2-6}
&$ 8\times 8 $	        &  4.2371e-04    &   2.97    & 1.9178e-05   &   2.96
\\ \cline{2-6}
&$16\times 16$	      &   5.3390e-05   &   2.99    &   2.4340e-06&   2.98
\\ \cline{2-6}
&$32\times 32$	      &  6.6913e-06    &  3.00     &3.5426e-07   &  2.98  
\\	\hline
\end{tabular}
\end{table}

\begin{table}[htbp]
\centering
\caption{History of convergence for Example \ref{Exa5.2} with $M=64$: temporal accuracy }
\label{table6}
\setlength{\tabcolsep}{4mm}{
\begin{tabular}{|c|c|c|c|c|c|c|c|}
	\hline
\multirow{8}{*}{$k=1$} 	&$\Delta t$    & $\frac{\lVert v(T)-v_{h0}(T) \rVert_0}{\lVert v(T)\rVert_0}$ & order &$\frac{\sqrt{\Delta t}||\varepsilon(v(T))-\varepsilon_h(v_h(T))||_0}{||\varepsilon(v(T))||_0}$&order 

\\ 	\hline
&$ 0.5 $   &  3.5128e-01    &--   &2.4872e-01 & --          

\\   \cline{2-6}
&$ 0.25 $   &  1.4792e-01    &1.25   &7.4196e-02 &1.75          

\\	 \cline{2-6}
&$ 0.125 $	&   6.7961e-02   &1.12     &2.4194e-02& 1.62  

\\	 \cline{2-6}
&$0.0625$	&  3.2583e-02    &1.06      &8.2635e-03 &  1.55
			
\\	 \cline{2-6}
&$0.03125$	&  1.5954e-02   &1.03    &2.9041e-03 &    1.51
			\\
			 \cline{2-6}
&			$0.015625$	&  7.8944e-03  &1.02     &1.0467e-03 &   1.47
			\\
			\hline
\multirow{6}{*}{$k=2$} 	&$ 0.5 $   &  3.5128e-01    &--   &2.4839e-01 & --          
\\
 \cline{2-6}

&$ 0.25 $   &  1.4792e-01    &1.25   &7.3962e-02 &1.75          
\\
 \cline{2-6}
&$ 0.125 $	&   6.7961e-02   &1.12     &2.4028e-02& 1.62  
\\
 \cline{2-6}
&$0.0625$	&  3.2583e-02    &1.06      &8.1460e-03 &  1.56
\\
 \cline{2-6}
&$0.03125$	&  1.5954e-02   &1.03    &2.8205e-03 &    1.53
\\
 \cline{2-6}
&$0.015625$	&  7.8944e-03  &1.02     &9.8680e-04 &   1.52
	\\
\hline			
			
	\end{tabular}}
\end{table}	

 \section{Conclusion}
In this paper, we have proposed a class of arbitrary order semi-discrete and fully-discrete WG finite element methods for the quasistatic Maxwell viscoelastic model. We have shown theoretically and numerically that the methods are of  optimal convergence rates. 

\bibliographystyle{plain}
\bibliography{mybibfile.bib}

\end{document}